\documentclass[a4paper,11pt]{article}
\author{Alan Hammond\\Department of Statistics\\University of California, Berkeley\\alanmh@stat.berkeley.edu}
\title{A lattice animal approach to percolation}

\newtheorem {theorem}{Theorem}[section]

\newtheorem {lemma}[theorem]{Lemma}

\newcounter{conjecture}
\newcounter{remark}
\newenvironment{conjecture}{\medskip{\bf Conjecture.\ }\em}{\rm}

\newcommand{\eqnsection}{
   \renewcommand{\theequation}{\thesection.\arabic{equation}}
   \makeatletter
   \csname @addtoreset\endcsname{equation}{section}
   \makeatother}

\newtheorem {definition}{Definition}[section]

\usepackage{amsfonts}
\usepackage[dvips]{graphicx}

\newcommand{\reta}{\varsigma}
\newcommand{\rrho}{\varrho}
\newcommand{\appak}{\lambda}
\newcommand{\lapha}{\alpha}
\newcommand{\sdash}{\sigma_{n,m}}
\newcommand{\sigm}{\sigma_{n,m}p^n(1-p)^m}
\newcommand{\signt}{\sigma_{n,m}t_n^n(1-t_n)^m}
\newcommand{\sigt}{{\Big( \frac{n}{n+m} \Big)}^{n} {\Big( \frac{m}{n+m} \Big)}^{m}}
\newcommand{\sigo}{p_c^n(1 - p_c)^m}

\newcommand{\suf}[2]{\sum_{n}{\sum_{#1}{#2}}}
\newcommand{\suh}[3]{\sum_{n}{\sum_{#1}^{#2}{#3}}}

\newcommand{\atwo}{m > \lfloor n \lapha + n^{1/2} \rfloor + 1}
\newcommand{\alphigh}{m \in n B ( \lapha , C ( \log{n} / n )^{1/2} ) }
\newcommand{\alplow}{m \in n B ( \lapha , n^{-1/2} ) }
\newcommand{\alpnhigh}{m \in n B ( \lapha_n , C ( \log{n} / n )^{1/2} ) }
\newcommand{\alpnlow}{m \in n B ( \lapha_n , n^{-1/2} ) }

\begin{document}
\maketitle
{\bf Abstract.}
We examine the percolation model on $\mathbb{Z}^d$ by an approach
involving lattice animals and their surface-area-to-volume ratio.
For $\beta \in [0,2(d-1))$, let
$f(\beta)$ be the asymptotic exponential rate in the number
of edges of the number of lattice animals containing the origin which
have surface-area-to-volume ratio $\beta$. The function $f$ is bounded
above by a function which may be written in an explicit form. For low
values of $\beta$ (\mbox{$\beta \leq 1/p_c - 1$}), equality holds, as
originally demonstrated by F.Delyon. For higher values ($\beta > 1/p_c
- 1$), the inequality is strict. 

We introduce two critical exponents, one of which describes how
quickly $f$ falls away from the explicit form as $\beta$ rises from
$1/p_c - 1$, and the second of which describes how large clusters
appear in the marginally subcritical regime of the percolation model.
We demonstrate that the pair of exponents must satisfy certain
inequalities, while other such inequalities yield sufficient
conditions for the absence of an infinite cluster at the critical
value. The first exponent is related to one of a more conventional nature in
the scaling theory of percolation, that of correlation size. In
deriving this relation, we find that there are two possible
behaviours, depending on the value of the first exponent, for the
typical surface-area-to-volume ratio of an unusually large cluster in
the marginally subcritical regime.

\newpage
\tableofcontents

\section{Introduction}

Percolation on the integer lattice $\mathbb{Z}^d$ is one of the most
fundamental and intensively studied models in the rigorous theory of
statistical mechanics. Many aspects of the behaviour of the model in
the subcritical and supercritical regime have been determined
rigorously. The problem of understanding the behaviour of the model at
criticality, and interplay between this behaviour and that for
parameter values nearby, has been addressed widely by physicists, but
the search for proofs of many of their predictions continues. These
predictions typically take the form of asserting the value of critical
exponents, and thereby describe the power-law decay or explosion of
characteristics of the model near criticality.   

In this report, we examine the percolation model by an
approach involving lattice animals, divided according to their
surface-area-to-volume ratio. Full details of proofs are provided, as
well as a discussion of some relevant literature. In a shortened
version \cite{pla} of this document, the central elements of
the approach and the proofs of the main theorems are presented. 
We will work throughout with the bond percolation
model in $\mathbb{Z}^d$. However, the results apply to the site or
bond model on any infinite transitive amenable graph with inessential changes.

For any given $p \in (0,1)$, two lattice animals with given size are
equally likely to arise as the cluster $C(0)$ containing the origin
provided that they have the same surface-area-to-volume ratio.
For given $\beta \in (0,\infty)$, there is an exponential growth rate
in the number of edges  for the number of lattice animals up to translation that have surface-area-to-volume ratio very close to $\beta$. This growth rate $f(\beta)$ may be studied as a function of $\beta$.
To illustrate the connection between the percolation model and the combinatorial question of the behaviour of $f$, note that the probability that the cluster containing the origin contains a large number $n$ of edges is given by
\begin{displaymath}
\mathbb{P}_p(\vert C(0) \vert = n) = \sum_{m}{\sigma_{n,m}p^n(1-p)^m},
\end{displaymath}
where $\sigma_{n,m}$ is the number of lattice animals that contain the origin, have $n$ edges and $m$ outlying edges. We rewrite the right-hand-side to highlight the role of the surface-area-to-volume ratio, $m/n$:
\begin{equation}\label{rrr}
\mathbb{P}_p(\vert C(0) \vert = n) = \sum_{m}{(f_n(m/n)p (1-p)^{m/n})^{n}}.
\end{equation}
Here $f_n(\beta) = (\sigma_{n,\lfloor \beta n \rfloor})^{1/n}$ is a
rescaling that anticipates the exponential growth that occurs. We
examine thoroughly the link between percolation and combinatorics
provided by (\ref{rrr}). 

An overview of the approach is now given, in the form of a
description of the organisation of the paper.
In Section \ref{secttwo}, we describe the model, and define notations, before
stating the combinatorial results that we will use. Theorem \ref{thmone}
asserts the existence and log-concavity of the function $f$, its
cumbersome proof appearing in the Appendix. Theorem \ref{thmtwo} implies
that
\begin{equation}\label{fgh}
\log f(\beta) \leq (\beta + 1)\log(\beta + 1) - \beta \log \beta \
\textrm{for $\beta \in (0,2(d-1))$}.
\end{equation}
F.Delyon \cite{DelyonT} showed that equality holds for $\beta \in (0,1/p_c -
1)$. Theorem \ref{thmtwo} implies that the inequality is strict for higher values of
$\beta$. 
 The marked change, as $\beta$ passes through $1/p_c - 1$, in the structure of large
lattice animals of surface-area-to-volume ratio $\beta$ is a
combinatorial analogue of the phase transition in percolation at criticality. 
The notion of a collapse transition for animals has been explored in \cite{MR1304212}.

In Section \ref{sectthr}, two scaling hypotheses are
introduced, each postulating the existence of a critical exponent. One of the
exponents, $\reta$, describes how quickly $f(\beta)$ drops away from
the explicit form given on the right-hand-side of (\ref{fgh}) as $\beta$ rises above $1/p_c - 1$. The other,
$\appak$, describes how rapidly decaying in $n$ is the discrepancy between the
critical value and that value on the subcritical interval at which the
probability of observing an $n$-edged animal as the cluster to
which the origin belongs is maximal. The first main result, Theorem
\ref{thmthr}, is then proved: the inequalities $\appak < 1/2$ and $\reta
\appak < 1$ cannot both be satisfied, because they imply that the mean
cluster size is uniformly bounded on the subcritical interval,
contradicting known results.

In Section \ref{sectthree}, sufficient conditions for the absence of an
infinite cluster at the critical value are proved. Theorem
\ref{thmthr} asserts that $\reta < 2$ or $\appak > 1/2$ are two such
conditions.
Except for some borderline cases, the range of values remaining after
Theorems \ref{thmthr} and \ref{thmfour} is specified by $\appak < 1/2$ and $\reta \appak >
1$. In Theorem \ref{thmund}, where we see that in this case, such a sufficient
condition may be expressed in terms of the extent to which the asymptotic
exponential rate $f(\beta)$ is underestimated by its finite
approximants $f_n(\beta)$ for a certain range of values of $\beta$.
The extent of underestimation is related to combinatorial exponents
such as the entropic exponent (see for example \cite{MR95m:82076}).

In Section \ref{sectfour}, we relate the value of $\reta$ to an
exponent of a more conventional nature in the scaling theory of
percolation, that of
correlation size (see Theorem \ref{thmfive}). 
Suppose that we perform an experiment in which the
surface-area-to-volume ratio of the cluster to which the origin
belongs is observed, conditional on its having a very large number of
edges, for a $p$-value slightly below $p_c$. 
How does the typical measurement, $\beta_p$, in this experiment behave 
as $p$ tends to $p_c$? 
The value $\beta_p$ tends to lie somewhere on the interval
\mbox{$(1/p_c - 1,1/p - 1)$}.
In Theorem \ref{thmsix}, we determine that there are two possible scaling
behaviours. The inequality $\reta < 2$ again arises, distinguishing the
two possibilities. If $\reta < 2$, then $\beta_p$ scales much closer to
$1/p_c - 1$ while if $\reta > 2$, it is found
to be closer to $1/p - 1$.

\section{Notations and combinatorial results}\label{secttwo}
\subsection{General definitions}
Throughout, we work with the bond percolation model on $\mathbb{Z}^d$,
for any given $d \geq 2$.
This model has a parameter $p$ lying in the interval $[0,1]$. Nearest
neighbour edges of $\mathbb{Z}^d$ are declared to be open with
probability $p$, these choices being made independently between
distinct edges. For any vertex $x \in \mathbb{Z}^d$, there is a
cluster $C(x)$ of edges accessible from $x$, namely the collection of
edges that lie in a nearest-neighbour path of open edges one of whose members contains
$x$ as an endpoint. The percolation probability $\theta(p)$ as a
function of $p$ may then be written \mbox{$\theta(p)= \mathbb{P}(
\vert C(0) \vert = \infty)$}. To demonstrate the continuity of
$\theta$, it suffices to show that $\theta(p_c) = 0$ (cf \cite{grim}), where $p_c$
denotes the critical value, namely the infimum of those values of $p$
for which $\theta$ is positive.

\begin{definition}
A lattice animal is the collection of edges of a finite connected
subgraph of $\mathbb{Z}^d$. An edge of $\mathbb{Z}^d$ is said to be
outlying to a lattice animal if it is not a member of the animal, and
if there is an edge in the animal sharing an endpoint with this
edge.
We adopt the notations:
\begin{itemize}
\item for $n,m \in \mathbb{N}$, set $\Gamma_{n,m}$ equal to the collection
of lattice animals in $\mathbb{Z}^d$ one of whose edges contains the
origin, having $n$ edges, and $m$ outlying edges. Define
$\sigma_{n,m} = \vert \Gamma_{n,m} \vert$. The surface-area-to-volume
ratio of any animal in $\Gamma_{n,m}$ is said to be $m/n$.
\item let $\Gamma'_{n,m}$ denote the subset of  $\Gamma_{n,m}$ whose
members' lexicographically minimal vertex is the origin. 
Set $\sigma'_{n,m} = \vert \Gamma'_{n,m} \vert$.
\item
for each $n \in \mathbb{N}$, define the functions $f_n,f_n': [0,\infty) \to [0,\infty)$ by
\begin{displaymath}
 f_n(\beta)=(\sigma_{n,\lfloor \beta n \rfloor})^{1/n},
 f'_n(\beta)=(\sigma'_{n,\lfloor \beta n \rfloor})^{1/n}
\end{displaymath} 
\end{itemize}
\end{definition}
On another point of notation, we will sometimes write the index set of
a sum in the form $n S$, with $S \subseteq (0,\infty)$, by which is meant
$\{m \in \mathbb{N}: m/n \in S \}$.
\subsection{Statement of combinatiorial theorems}
We require some results about the asymptotic exponential growth rate
 of the number of lattice animals as a function of their
 surface-area-to-volume ratio.
\begin{theorem} \label{thmone} 
${}$
\flushleft
\begin{enumerate}
\item For $\beta \in [0,\infty) - \{ 2(d-1) \}$, $f(\beta)$ exists,
being defined as the limit $\lim_{n \to \infty} f_n(\beta) $. 
\item   for $ \beta > 2(d-1)$, $f(\beta) = 0$. 
\item for $\beta \in (0,2(d-1)), n \in
\mathbb{N}$, $f_n$ satisfies $f_n(\beta) \leq L^{1/n} n^{1/n}
 f(\beta)$, where the constant $L$ may be chosen uniformly in $\beta \in (0,2(d-1))$.
\item f is log-concave on the interval $(0,2(d-1))$.
\end{enumerate}
\end{theorem}
\begin{theorem}\label{thmtwo}
Introducing $ g: (0,2(d-1)) \to [0,\infty)$ by means of the formula
\begin{displaymath}
 f(\beta) = g(\beta)\frac{ (\beta +1)^{\beta +1} } { \beta^\beta },
\end{displaymath}
we have that
\begin{displaymath}
g(\beta) \left\{ \begin{array}{ll}
        = 1 & \textrm{on $(0,\lapha]$,}\\
        < 1 & \textrm{on $(\lapha,2(d-1))$,}
\end{array} \right. 
\end{displaymath} 
where throughout $\lapha$ denotes the value $ 1/p_c - 1$.       
\end{theorem}
{\bf Remark} The assertion that $g = 1$ on $(0,\lapha]$ was originally
proved by Delyon \cite{DelyonT}.

The proof of Theorem \ref{thmone} of $f$ is a lengthy and tiresome
task, which is performed in the Appendix.
\subsection{Proof of Theorem $\ref{thmtwo}$}
\subsubsection{Some lemmas}
The following lemmas will be used.
\begin{lemma}\label{lemmasix}
$f: (0,2(d-1)) \to [0,\infty)$ satisfies
\begin{displaymath} 
\log f(\beta) \leq (\beta +1) \log (\beta +1) - \beta \log \beta
\end{displaymath}
\end{lemma} 
{\bf Proof.}
We give a probabilistic proof, that uses the percolation model.
The probability that the cluster to which the origin belongs contains $n$ edges can be represented as a sum of terms including the expressions $\sigma_{n,m}$:
\begin{displaymath}
\mathbb{P}_p (\vert C(0) \vert = n) = \sum_{m}{\sigma_{n,m}p^{n}(1-p)^m}
\end{displaymath}
Let $\beta \in (0,2(d-1))$. Choosing $p = 1/(1 + \beta)$, and noting that the right-hand-side of the above equation is bounded above by one, yields
\begin{displaymath}
\sigma_{n,\lfloor \beta n \rfloor} \leq {\bigg( \frac{(\beta + 1)^{\beta +1}}{\beta^{\beta}}\bigg) }^n ,
\end{displaymath}
which may be rewritten
\begin{displaymath}
f_n(\beta) \leq \frac{(\beta + 1)^{\beta +1}}{\beta^{\beta}}.
\end{displaymath}
Taking the limit as $n \to \infty$ gives that
\begin{displaymath}
f(\beta) \leq \frac{(\beta + 1)^{\beta +1}}{\beta^{\beta}},
\end{displaymath}
as required. $\Box$ 
\begin{lemma}\label{lemmanex}
There exists $r \in (0,1)$ such that, for $n$ sufficiently large and for  $m \in \{ 2(d-1)n, \ldots, 2(d-1)n + 2d \}$, we have that
\begin{displaymath} 
\sigma_{n,m} \leq  \frac{ ( 1+ \frac{m}{n} )^{n+ m}}{(\frac{m}{n})^{m}} r^n . 
\end{displaymath}
\end{lemma} 
{\bf Proof.} Let $j \in \{ 0, \ldots, 2d \}$. Then
\begin{eqnarray}
  & & \sigma_{n,2(d-1)n + j} \Big( \frac{n}{(2d-1)n + j} \Big)^n \Big( \frac{2(d-1)n + j}{(2d-1)n + j} \Big)^{2(d-1)n + j} \nonumber \\
 & & \quad \leq \, \mathbb{P}_{n/{((2d-1)n + j)}}{(\vert C(0) \vert = n)} . \nonumber
\end{eqnarray}
The inequality $p_c > 1/(2d-1)$ is of course known, though in an aside
to the notes for this section, a proof is supplied. There is an exponential decay rate in $n$ for the probability of observing an $n$-edged cluster containing the origin in the subcritical phase \cite{MR86h:82045}. Hence, for some $r_j \in (0,1)$ and for $n$ sufficiently large, we have that
\begin{displaymath}
\sigma_{n,2(d-1)n + j} \leq \Big( \frac{(2d-1)n + j}{n} \Big)^n \Big( \frac{(2d-1)n + j}{2(d-1)n + j} \Big)^{2(d-1)n + j} r_j^n .
\end{displaymath}
Setting $r = \max_{j \in \{0, \ldots, 2d \}}{r_j}$ gives the result. $\Box$ \\ 
\subsubsection{Proof of Delyon's result.} 
For completeness, we include a proof
of Delyon's result.
We know that $f$ is log-concave on $(0,\lapha)$ and that on that interval, it satisfies 
\begin{displaymath}
f(\beta) \leq  \frac{(\beta + 1)^{\beta +1}}{\beta^{\beta}} .
\end{displaymath}
From these statements and the assumption that Delyon's result fails, it follows
that there exists $\beta_0 \in (0,\lapha), \epsilon > 0$ such that 
\begin{displaymath}
f(\beta) \leq \frac{(\beta + 1)^{\beta +1}}{\beta^{\beta}} - \epsilon \quad \textrm{on} \, \, (\beta_0 - \epsilon,\beta_0 + \epsilon).
\end{displaymath}
Set $p = 1/(1 + \beta_0)$, and note that $p > p_c$.
We have that
\begin{displaymath}
\mathbb{P}_p (\vert C(0) \vert = n) = \sum_{m}{\sigma_{n,m}p^{n}(1-p)^m}.
\end{displaymath}
Note that
\begin{eqnarray}
& & \sum_{m}{\sigma_{n,m}p^{n}(1-p)^m} \nonumber \\
& = & \sum_{m}{{\Big( f_n(m/n)p(1-p)^{\frac{m}{n}} \Big) }^n} \nonumber \\
& = & \sum_{m}{{\bigg( f_n(m/n)  \frac{ {\beta_0}^{\frac{m}{n}}}{{(1 + \beta_0)}^{1 + \frac{m}{n}}} \bigg)}^n}. \nonumber
\end{eqnarray}
Then 
\begin{eqnarray}
\mathbb{P}_p (\vert C(0) \vert = n) & = & \sum_{m \in n S_1}{{\bigg( f_n(m/n)  \frac{ {\beta_0}^{\frac{m}{n}}}{{(1 + \beta_0)}^{1 + \frac{m}{n}}} \bigg)}^n} {} \nonumber \\
& & + \sum_{m \in n S_2}{{\bigg( f_n(m/n)  \frac{ {\beta_0}^{\frac{m}{n}}}{{(1 + \beta_0)}^{1 + \frac{m}{n}}} \bigg)}^n} \nonumber \\
& &  + \sum_{m = 2(d-1)n}^{2(d-1)n + 2d}{{\bigg( f_n(m/n)  \frac{ {\beta_0}^{\frac{m}{n}}}{{(1 + \beta_0)}^{1 + \frac{m}{n}}} \bigg)}^n} \nonumber
\end{eqnarray}
where 
\begin{displaymath}
S_1 = (\beta_0 - \epsilon,\beta_0 + \epsilon)  \quad \textrm{ and } \quad
S_2 = \big( 0,2(d-1) \big) - (\beta_0 - \epsilon,\beta_0 + \epsilon).
\end{displaymath}
The behaviour of the three sums above will now be analysed, under the
assumption that Delyon's result is false. Firstly, we need a definition.
\begin{definition}\label{defnphi}
Let the function $\phi: [0,\infty)^2 \to \mathbb{R}$ be given by
\begin{equation}
\phi \big( \beta,\gamma \big) = (\gamma + 1)\log(\gamma + 1) - \gamma \log \gamma + \gamma \log \beta - (\gamma + 1)\log(\beta + 1).
\end{equation}
\end{definition}
\begin{itemize}
\item{The sum indexed by $n S_1$}

We know from Theorem $\ref{thmone}$ (iii) that there exists $L > 1$
such that, for all $n \in \mathbb{N}$, $\beta \in (0,2(d-1))$,
$f_n(\beta) \leq (Ln)^{1/n} f(\beta)$. It follows that
\begin{eqnarray}
& & \sum_{m \in n S_1}{{\bigg( f_n(m/n)  \frac{ {\beta_0}^{\frac{m}{n}}}{{(1 + \beta_0)}^{1 + \frac{m}{n}}}   \bigg)}^n} \nonumber \\
& \leq & L n \sum_{m \in n S_1}{{\bigg( f(m/n)  \frac{ {\beta_0}^{\frac{m}{n}}}{{(1 + \beta_0)}^{1 + \frac{m}{n}}}   \bigg)}^n} \nonumber \\
& \leq &  L n  \sum_{m \in n S_1}{{ \Bigg( \bigg( \frac{ ( 1+ \frac{m}{n} )^{1+\frac{m}{n}}}{(\frac{m}{n})^{\frac{m}{n}}} -\epsilon  \bigg) \bigg(  \frac{ {\beta_0}^{\frac{m}{n}}}{{(1 + \beta_0)}^{1 + \frac{m}{n}}} \bigg) \Bigg)}^n } \nonumber \\
& \leq & L n \sum_{m \in n S_1}{\exp{n[\phi(\beta_0,m/n) + \log(1 - \epsilon R)]}} \nonumber \\
& \leq & L n (2 \epsilon n +1)(1- \epsilon R)^n, \nonumber
\end{eqnarray}
where $R$ is some positive constant.

\item{The sum indexed by $n S_2$}

In this case, note that
\begin{eqnarray}
& & \sum_{m \in n S_2}{{\bigg( f_n(m/n)  \frac{ {\beta_0}^{\frac{m}{n}}}{{(1 + \beta_0)}^{1 + \frac{m}{n}}}   \bigg)}^n} \nonumber \\
& \leq & L n \sum_{m \in n S_2}{{\bigg( f(m/n)  \frac{ {\beta_0}^{\frac{m}{n}}}{{(1 + \beta_0)}^{1 + \frac{m}{n}}}   \bigg)}^n} \nonumber \\
& \leq &  L n \sum_{m \in n S_2}{{ \Bigg( \bigg( \frac{ ( 1+ \frac{m}{n} )^{1+\frac{m}{n}}}{(\frac{m}{n})^{\frac{m}{n}}} \bigg) \bigg(  \frac{{\beta_0}^{\frac{m}{n}}}{(1+ \beta_0)^{1 + \frac{m}{n}}} \bigg) \Bigg)  }^n } \nonumber \\
& = & L n \sum_{m \in n S_2}{\exp{n \phi (\beta_0,m/n) }}. \nonumber 
\end{eqnarray}
The fact that
\begin{displaymath}
 \frac{d}{d \gamma} \phi (\beta_0,\gamma) =  \log (1+ 1/{\gamma}) - \log(1 + 1/{\beta_0})
\end{displaymath}
implies that there exists $\delta > 0$ such that $\phi(\beta_0,\gamma)
< - \delta$ for $\gamma \in S_2$.
Hence
\begin{displaymath}
\sum_{m \in n S_2}{{\bigg( f_n(m/n)  \frac{
{\beta_0}^{\frac{m}{n}}}{{(1 + \beta_0)}^{1 + \frac{m}{n}}}
\bigg)}^n} \leq 2(d-1) L n^2 \exp{-n\delta}
\end{displaymath}

\item{The sum indexed by $\{2(d-1)n, \ldots, 2(d-1)n + 2d \}$}

Note that Lemma \ref{lemmanex} implies that, for $n$ sufficiently large,
\begin{eqnarray}
& & \sum_{m = 2(d-1)n}^{2(d-1)n + 2d}{{\bigg( f_n(m/n)  \frac{ {\beta_0}^{\frac{m}{n}}}{{(1 + \beta_0)}^{1 + \frac{m}{n}}}   \bigg)}^n} \nonumber \\
& \leq & \sum_{m = 2(d-1)n}^{2(d-1)n + 2d}{{ \Bigg( \bigg( \frac{ ( 1+ \frac{m}{n} )^{1+\frac{m}{n}}}{(\frac{m}{n})^{\frac{m}{n}}} \bigg) \bigg(  \frac{{\beta_0}^{\frac{m}{n}}}{(1+ \beta_0)^{1 + \frac{m}{n}}} \bigg) \Bigg)  }^n } r^n \nonumber \\
& = & \sum_{m = 2(d-1)n}^{2(d-1)n + 2d}{r^n \exp{n\phi (\beta_0,m/n)} }. \nonumber \\
& \leq & (2d+1) r^n . \nonumber  
\end{eqnarray}
\end{itemize}
We have demonstrated if Delyon's result fails, then 
\begin{displaymath}
\liminf_{n \to \infty}{\frac{- \log \mathbb{P}_p(\vert C(0) \vert = n)}{n}} > 0 .
\end{displaymath}

The subexponential decay rate for the probability of observing a large cluster in the supercritical phase was proved by H. Kunz and B. Souillard \cite{MR58:14852} and by M. Aizenman, F. Delyon and B. Souillard \cite{MR82b:82048}.
Since $p > p_c $ is in the supercritical phase, the above expression is zero.
This contradiction completes the proof of Delyon's result.

\subsubsection{$g<1$ for $\beta > \lapha$}

We aim to show that, for $\beta \in (\lapha,2(d-1))$, $g(\beta)$ is strictly less than one.
Let $\beta$ lie in this interval. Let $p = 1/(1 + \beta)$. Note that $p < p_c$, and that
\begin{eqnarray}
\mathbb{P}_p(\vert C(0) \vert = n) & \geq & \mathbb{P}_p(C(0) \in \Gamma_{n,\lfloor \beta n \rfloor}) \nonumber \\
 & = & \vert \Gamma_{n,\lfloor \beta n \rfloor} \vert \frac{{\beta}^{\lfloor \beta n \rfloor}}{(1 + \beta)^{n + \lfloor \beta n \rfloor}} \nonumber \\
 & = &  \big( f_n(\beta) \big)^n \frac{{\beta}^{\lfloor \beta n \rfloor}}{(1 + \beta)^{n + \lfloor \beta n \rfloor }}. \nonumber 
\end{eqnarray}
Taking logarithms yields
\begin{displaymath}
 \frac{\log \mathbb{P}_p \big( \vert C(0) \vert = n \big)}{n} \geq
 \log f_n(\beta) + \frac{\lfloor \beta n \rfloor \log \beta}{n} -
 \Big( 1 + \frac{\lfloor \beta n\rfloor }{n} \Big) \log (1 + \beta) ,
\end{displaymath}
from which it follows that
\begin{equation}\label{pop}
 \liminf_{n \to \infty}{\frac{\log \mathbb{P}_p \big( \vert C(0) \vert = n \big)}{n}} \geq  \log f(\beta) + \beta \log \beta - (1+ \beta) \log (1+ \beta) .
\end{equation}
The right-hand-side of (\ref{pop}) is equal to $\log g(\beta)$, by definition.
The exponential decay rate for the probability of observing a large cluster in the subcritical phase was established by M. Aizenman and C.M. Newman in \cite{MR86h:82045}. Since $p < p_c$, this means the left-hand-side of (\ref{pop}) is negative. This implies that $g(\beta) < 1$, as required.
This completes the proof of Theorem $\ref{thmtwo}$. $\Box$

\subsection{Notes on properties of lattice animals}
We recount the work of various authors on properties of lattice animals in the integer lattice. 

The form of the function $f$ on the interval $(0,\lapha)$ was
established in Theorem \ref{thmtwo}. This result was discovered by
F. Delyon in his thesis \cite{DelyonT}. Theorem $II.2$ on page $24$ of
\cite{DelyonT} asserts: \\
{\bf Theorem.}
Let $m$ and $n$ satisfying $0 < m/n = k < (1 - p_c)/p_c$ be given. Then,
\begin{displaymath}
\log \frac{(1+k)^{1+k}}{k^k} \geq \frac{\log \sigma_{n,m}}{n} \geq \log \frac{(1+k)^{1+k}}{k^k} + O(n^{-1/{2 d}}).
\end{displaymath}

In \cite{MR90c:82051}, Madras et al. investigate the large-$n$ behaviour of the number of lattice animals with $n$ vertices and $\lapha n$ cycles. The authors define $a_n(c)$ to be the number (up to translation) of lattice animals with $n$ vertices and cyclomatic index $c$. That is, $c$ is the minimal number of edges whose removal produces a tree.
For $\lapha \in [0,d-1)$, the existence of the limit
\begin{displaymath}
 \phi(\lapha) = \lim_{n \to \infty}{a_n(\lceil n \lapha \rceil)^{1/n}},
\end{displaymath}
 is demonstrated. As in the proof of Theorem \ref{thmone}, relegated
 to the Appendix, a pair of large lattice animals is concatenated, and a slight discrepancy in the number of cycles of the new animal occurs. The authors are able to dispense with the protracted considerations of a correction construction by demonstrating that,

\begin{displaymath}
{{c + k} \choose k} a_n(c+k) \leq {{(d-1)n - c - 1} \choose k} a_n(c) \ \textrm{for all $n$}.
\end{displaymath} 

This result demonstrates the subadditivity of $a_n(\lceil n\lapha \rceil)^{1/n}$ and obviates the need for some other construction to effect the slight adjustments to the parameter $c$ in the animal formed in the course of the concatenation argument.

After making the definitions 
\begin{displaymath}
A_n(z)  =  \sum_{c}{a_n(c)z^c}
\end{displaymath}
and
\begin{displaymath}
\log \Lambda (z)  =  \lim_{n \to \infty}{\frac{\log A_n (z)}{n}},
\end{displaymath}
the authors prove that
\begin{equation}\label{eqnadh}
\log \Lambda (z) = \sup_{0 \leq \lapha \leq d-1}{\log \phi(\lapha) + \lapha \log z}.
\end{equation} 
That is, if an experiment is conducted in which an $n$-edged lattice animal is chosen according to a distribution which penalises each additional cycle by a factor of $z$ ( $= 1 - p$, for $p \in (0,1)$, say), then the likely number of cycles of the observed animal will be of the form $\lapha n$, where $\lapha$ attains the supremum in (\ref{eqnadh}).

In \cite{MR2001f:82034}, N.Madras addresses the following question. Let $B$ be small integer. Let $\chi$ denote any configuration of edges lying in the box $\tau_B = \{ -B, \ldots, B \}^d$ such that any edge of $\chi$ may be reached from the boundary of $\tau_B$ by a path of edges of $\chi$. 
For $\gamma$ a lattice animal with a large number $n$ of edges, let $\gamma_{\chi}$ be equal to the proportion of vertices $v$ of $\gamma$ such that $\{ \tau_B + v \} \cap \gamma$ is equal to the given pattern $\chi$. The quantity $\gamma_{\chi}$ should be thought of as the proportion of instances of the pattern $\chi$ in the large animal $\gamma$. Madras proves that there exists $\epsilon > 0$ such that,
\begin{displaymath} 
\limsup_{n \to \infty}{{\vert \{ \gamma \in \Gamma_n : \gamma_{\chi} < \epsilon \} \vert}^{1/n}} < \lim_{n \to \infty}{{\vert \Gamma_n \vert}^{1/n}},
\end{displaymath}
for any such pattern of edges $\chi$ lying in the box $\tau_B$.
The statistic $\gamma_{\chi}$ performs the role that surface-area-to-volume ratio does in this report. Madras' result is a step on the road to proving existence, and  the strict log-concavity of the function
\begin{displaymath}
 [0,1] \to \mathbb{R}^{+}: x \to \lim_{n \to \infty}{{\vert \{ \gamma
 \in \Gamma_n : \gamma_{\chi} = \lfloor n x \rfloor  \} \vert}^{1/n}}.
\end{displaymath}
The analogous question for the approach being developed in this report is the strict log-concavity of $f$. This question appears to be unsettled. 

Flesia et al \cite{MR1304212} study a two-variable model in which the outlying edges of a lattice animal are divided into two camps. There are the contacts, which are those outlying edges both of whose endpoints lie in the animal, and the solvents, which are the remainder. Writing $a_n(s,k)$ for the number (up to translation) of animals with $n$ vertices, $s$ solvents and $k$ contacts, the two variable partition function
\begin{displaymath}
Z_n(\beta_1,\beta_2) = \sum_{s,k}{a_n(s,k)\exp{(s \beta_1 + k \beta_2)}},
\end{displaymath}
is defined, and the limiting free-energy, given by,
\begin{displaymath}
G(\beta_1,\beta_2) = \lim_{n \to \infty}{\frac{Z_n(\beta_1,\beta_2)}{n}},  
\end{displaymath}
is shown to exist, for all $(\beta_1,\beta_2) \in [0,\infty)^2$. The function $G$ is proved to be a convex and continuous function of its variables.
Van Rensberg et al \cite{ROT} consider the behaviour of the function
\begin{displaymath}
Z_n'(\beta_c,\beta_k) = \sum_{s,k}{a_n'(c,k)\exp{(c \beta_c + k \beta_k)}},
\end{displaymath}
where $a_n'(c,k)$ is the number (up to translation) of animals with $n$ vertices, $k$ contacts and cyclomatic index $c$.
The relation,
\begin{equation}\label{reln}
s + 2k + 2c = 2(d-1)n + 2,
\end{equation}
means that $Z_n'$ may easily be expressed in terms of $Z_n$. Relations like (\ref{reln}) allow the probability that the cluster containing the origin, $C(0)$, has $n$-edges in the percolation model with parameter $p$ to be expressed in the following form:
\begin{displaymath}
\mathbb{P}_p(\vert C(0) \vert = n) = p^n(1-p)^{(2d + 2(d-1)n)} \sum_{c,k}{(n+1 - c)a_n(c,k)(1-p)^{-(2dc + k)}}.
\end{displaymath}
In this form, $\mathbb{P}_p(\vert C(0) \vert = n)$ is a weighted sum of $Z_n'$ and its derivatives, provided that the choices
\begin{displaymath}
\beta_c = -2d \log (1-p); \, d \beta_k = - \log (1-p)
\end{displaymath}
are made. So a one-dimensional subset of the $(\beta_c,\beta_k)$-plane represents the percolation model, and this line may be indexed by the parameter $p$ of the model.
For fixed choices for $\beta_c$ and $\beta_k$, an experiment may be
performed in which an animal with $n$ vertices, $k$ contacts and
cyclomatic index $c$, is sampled with weight $\exp(c \beta_c + k
\beta_k)$. The quantity $Z_n'$ is then the partition function, that
is, the value by which we must normalize to obtain a probability measure.
The authors of \cite{ROT} explore the idea that a `collapse
transition' occurs as the parameters $\beta_c$ and $\beta_k$
vary. They hypothesise that for low values of the parameters, a large
animal sampled with this weighting has an expanded form, with, for example, comparatively high mean perimeter. This hypothesis is advanced by numerical evidence. The one-parameter subset corresponding to the percolation model for parameter values $p \in (0,p_c)$ is hypothesised to mark the boundary between two parts of the `collapsed' section of the $(\beta_c,\beta_k)$-parameter space. On one side, typical large animals are rich in cycles, and, on the other, they are rich in contact edges. 
For more work relating to the idea of collapse transition, see \cite{flesiaetal} and \cite{MR1421925}. 

In our context, the parameter that describes lattice animals is the surface-area-to-volume ratio. Lattice animals with low surface-area-to-volume ratio will typically occur for high values of the percolation parameter $p$; indeed, for $p > p_c$ lying in the supercritical phase, Delyon's result on the form of the function $f$ easily implies that large finite clusters in the model have surface-area-to-volume ratio close to $\beta = 1/p - 1$.
The collapse transition examined for two-variable models by \cite{MR1304212} and \cite{ROT} is also apparent in the approach adopted here. If surface-area-to-volume ratio $\beta$ takes a value lying in the interval $(\lapha,2(d-1))$, then a typical large animal with this parameter value looks like an unusually large subcritical percolation cluster, and has a relatively dispersed or stringy appearance. If $\beta$ takes the value $\lapha$, then the large animal looks like a critical percolation cluster. As $\beta$ falls below $\lapha$, the collapse transition occurs, and the large animal has the appearance of a ball whose interior resembles a region of the supercritical infinite cluster. As such, the animal is more compact than its higher-$\beta$ counterparts.
The form of the large finite cluster in the supercritical regime is
the subject of Cerf's monograph \cite{MR1774341}. \\
{\bf Aside: a proof of lower bound on $p_c$.}

In the proof of Lemma \ref{lemmanex}, we made use of the strict inequality
$p_c > 1/(2d-1)$. For completeness, here is a proof of this assertion.

\begin{lemma}\label{strpc}
For $d \geq 2$, the critical value $p_c$ satisfies the lower bound
\begin{displaymath}
p_c \geq \frac{1}{(2d-1) \big( 1 - [1/(2d - 1)]^2 \big)^{1/3}} .
\end{displaymath}
\end{lemma}   
{\bf Proof} Let $N_k$ denote the number of open paths containing $k$
edges that start at the origin. Since the event that $C(0)$ is infinite
implies that $N_k \geq 1$, we have that $\theta(p) \leq
\mathbb{E}_{p}(N_k)$. The expectation $\mathbb{E}_{p}(N_k)$ is bounded
above by $p^k R_k$, where $R_k$ is the number of paths that contain
$k$ edges starting from $0$. We now demonstrate that, for $k \in
\mathbb{N}$,
\begin{equation}\label{asser}
R_{3k+1} \leq 2d \Big[ (2d-1)^3 \big( 1 - (1/(2d-1))^2 \big) \Big]^k ,
\end{equation}
from which, the statement of the lemma follows immediately.
Let us prove (\ref{asser}) by an induction, from the trivial case of
$k=0$. Let a path $\tau$ with $3k + 1$ edges be given. The number of ways
of extending it to a path with three more edges is bounded above by
\begin{displaymath}
(2d)^3 - (2d)^2 - (2d-1)2d - (2d-1)^2 - (2d-1) . 
\end{displaymath}   
The successive negative contributions arise from the fact that the
path may not backtrack in its first added edge; it may not backtrack with
its second but not its first added edge; it may not backtrack with
only its third edge; the three added edges and the $3k+1$-st edge of
$\tau$ may not form a square.
This observation extends the induction and proves the assertion
(\ref{asser}). $\Box$

\section{Critical exponents and inequalities}\label{sectthr}

We introduce scaling hypotheses, which propose the existence of some
critical exponents. There are two scaling hypotheses we require, one
of which asserts the existence of an exponent describing the behaviour
of $f$ for values of the argument just greater than $\lapha$, and the
second of which describes the nature of the formation of large
clusters just below criticality. We then state and prove the first
main theorem, which demonstrates that a pair of inequalities involving
the two exponents cannot both be satisfied.

\subsection{Hypotheses on the existence of exponents}
\subsubsection{Hypothesis $(\appak)$}
\begin{definition}
Let $\{ t_n: n \in \mathbb{N} \}$ be a collection of values in $(0,p_c)$ satisfying the condition
\begin{equation}\label{eteen}
\sum_{m}{\sdash t_n^n (1 - t_n)^m} = \sup_{p \in (0,p_c]}{\sum_{m}{\sdash p^n (1 - p)^m}}. 
\end{equation}
\end{definition}
That is, $t_n$ is some point at or below the critical value at which the probability of observing an $n$-edged animal as the cluster to which the origin belongs is maximal. For definiteness, if there is more than one value with this property, we choose $t_n$ to be the least such. It is reasonable to suppose that $t_n$ is slightly less than $p_c$, and that the difference decays polynomially in $n$ as $n$ tends to infinity.

\begin{definition} 
Define $\Omega_{+}^{\appak} = \{ \beta \geq 0 : \liminf_{n \to \infty}{(p_c - t_n)/n^{-\beta}} = \infty \}$, and 
$\Omega_{-}^{\appak} = \{ \beta \geq 0 : \limsup_{n \to \infty}{(p_c - t_n)/n^{-\beta}} = 0 \}$. 
\end{definition}
\begin{definition}
If $\sup{\Omega_{-}^{\appak}} = \inf{\Omega_{+}^{\appak}}$, then hypothesis $(\appak)$ is said to hold, and $\appak$ is defined to be equal to the common value.
\end{definition}
So, if hypothesis $(\appak)$ holds, then $p_c - t_n$ behaves like $n^{-\appak}$, for large $n$. We remark that it would be consistent
with the notion of a scaling window about criticality that the probability of observing
the cluster $C(0)$ with $n$-edges achieves its maximum on the
subcritical interval on a short plateau whose right-hand endpoint is
the critical value. If this is the case, then $t_n$ should lie at the
left-hand endpoint of the plateau. To be confident that $p_c - t_n$ is
of the same order as the length of this plateau, the definition of the
quantities $t_n$ could be changed, so that a
small and fixed constant multiples the right-hand-side of  
(\ref{eteen}). In this report, any proof of a statement involving the
exponent $\appak$ is valid if it is defined in terms of this altered
version of the quantities $t_n$. 

\subsubsection{Hypothesis $(\reta)$}

This hypothesis is introduced to describe the behaviour of $f$ for values of the argument just greater than $\lapha$. It has been shown in Section \ref{secttwo} that the value $\lapha$ is the greatest for which $\log f(\beta) = (\beta +1)\log(\beta + 1)- \beta \log \beta$; the function $g$ was introduced to describe how $\log f$ falls away from this function as $\beta$ increases from $\lapha$. Thus, we phrase hypothesis $(\reta)$ in terms of $g$.
\begin{definition} 
Define $\Omega_{-}^{\reta} = \{ \beta \geq 0 : \liminf_{\delta \to 0}{(g(\lapha + \delta)- g(\lapha))/{\delta}^{\beta}} = 0 \}$, and 
$\Omega_{+}^{\reta} = \{ \beta \geq 0 : \limsup_{\delta \to 0}{(g(\lapha + \delta) - g(\lapha))/{\delta}^{\beta}} = - \infty \}$.
\end{definition}
\begin{definition}
If $\sup{\Omega_{-}^{\reta}} = \inf{\Omega_{+}^{\reta}}$, then hypothesis $(\reta)$ is said to hold, and $\reta$ is defined to be equal to the common value.
\end{definition}
{\bf Remark} It follows from Theorem \ref{thmtwo} that $g(\lapha)=1$.\\
If hypothesis $(\reta)$ holds, then greater values of $\reta$ correspond to a smoother behaviour of $f$ at $\lapha$. For example, if $\reta$ exceeds $N$ for $N \in \mathbb{N}$, then $f$ is $N$-times differentiable at $\lapha$.

\subsection{Exponent inequalities: the proof of Theorem \ref{thmthr}}
\begin{theorem}\label{thmthr}
Suppose that hypotheses ($\reta$) and ($\appak$) hold.  If $\appak < 1/2$, then $\reta \appak \geq 1$.
\end{theorem}
{\bf Proof}
We prove the Theorem by contradiction, assuming that the two hypotheses hold, and that $\appak < 1/2$, $ \reta \appak < 1$. Our aim is to arrive at the conclusion that the mean cluster size, given by $\sum_{n}{n \mathbb{P}_{p}{(\vert C(0) \vert = n)}}$, is bounded above, uniformly for $p \in (0,p_c)$. This contradicts a result of Aizenman and Newman (see \cite{MR86h:82045}). 
Note that
\begin{displaymath}
\sup_{p \in (0,p_c)}{\sum_{n}{n \mathbb{P}_{p}{(\vert C(0) \vert = n)}}} \leq \sum_{n}{n \mathbb{P}_{t_n}{(\vert C(0) \vert = n)}}.
\end{displaymath}
We write 
\begin{equation}\label{enst}
\mathbb{P}_{t_n}{(\vert C(0) \vert = n)} = \sum_{m}{\sigma_{n,m} t_n^n(1 - t_n)^m},
\end{equation}
and split the sum on the right-hand-side of (\ref{enst}). To do so, we use the following definition.
\begin{definition}\label{alphn}
For $n \in \mathbb{N}$, let $\lapha_n$ be given by $t_n = 1/(1 + \lapha_n)$.
For $G \in \mathbb{N}$, let  $D_n ( = D_n(G))$ denote the interval 
\begin{displaymath}
D_n = ( \lapha_n - G{\{\log(n)/n \}}^{1/2} , \lapha_n + G{\{\log(n)/n \}}^{1/2}).
\end{displaymath} 
\end{definition}
Now, 
\begin{displaymath} 
\sum_{m}{\sigma_{n,m} t_n^n(1 - t_n)^m} = C_1(n) + C_2(n) + C_3(n) ,
\end{displaymath}
where the terms on the right-hand-side are given by 
\begin{eqnarray}
 C_1(n) & = & \sum_{m \in n D_n}{\sigma_{n,m}t_n^n(1 - t_n)^m} , \nonumber \\
 C_2(n) & = &  \sum_{m \in n \big( (0,2(d-1)) - D_n
 \big)}{\sigma_{n,m}t_n^n(1 - t_n)^m}, \nonumber
\end{eqnarray}
and
\begin{displaymath}
 C_3(n)  =   \sum_{m \in \{ 2(d-1)n ,\ldots, 2(d-1)n + 2d \}}{\sigma_{n,m}t_n^n(1 - t_n)^m}. 
\end{displaymath}
\begin{lemma}\label{lemfour}
The function $\phi$ specified in  Definition \ref{defnphi} satisfies
\begin{displaymath}
 \phi \big( \lapha , \lapha + \gamma \big) = - \frac{\gamma^2}{2\lapha(\lapha + 1)} + O(\gamma^3).
\end{displaymath}
\end{lemma}
We compute
\begin{eqnarray}
\phi(\lapha,\lapha + \gamma) & = & - (\lapha + \gamma)\log(\lapha + \gamma) + (\lapha + 1 + \gamma)\log (\lapha + 1 + \gamma) {} \nonumber \\
   & & {} - (\lapha + 1 + \gamma)\log (\lapha + 1)  + (\lapha + \gamma) \log \lapha {} \nonumber \\
  & = & - (\lapha + \gamma)\log(1 + \gamma/{\lapha}) \nonumber \\
  & & {} + (\lapha + 1 + \gamma)\log (1 + \gamma/{(\lapha + 1)}) {} \nonumber \\
  & = & - (\lapha + \gamma)[\gamma/{\lapha} - \gamma^2/{2{\lapha}^2} ] \nonumber \\
  & & {} + (\lapha + 1 + \gamma)[\gamma/{(\lapha + 1)} - \gamma^2/{2(\lapha +1)^2} ]  + O(j(n)^3) {}\nonumber \\
  & = & - \gamma^2/[2\lapha(\lapha + 1)] + O(\gamma^3), {}\nonumber  
\end{eqnarray}
giving the result. $\Box$. \\
We have that
\begin{eqnarray}
 \sum_{n}{C_2(n)} & = & \suf{m \in n((\lapha,2(d-1)) - D_n)}{\sdash t_n^n(1 - t_n)^m} \nonumber \\
 & = & \suf{m \in n((\lapha,2(d-1)) - D_n)}{{\bigg( f_n(m/n)\frac{\lapha_n^{m/n}}{(1 + \lapha_n)^{1 + m/n}} \bigg)}^n} \nonumber \\
 & \leq & L \sum_{n}{ n  \sum_{m \in n((\lapha,2(d-1)) - D_n)}{\exp{\{n\phi_{\lapha_n,m/n}\}}}}, \nonumber 
\end{eqnarray}
where the inequality is valid by virtue of Theorem \ref{thmone} and
the fact that $g \leq 1$. Lemma \ref{lemfour} implies that
\begin{displaymath}
 \sum_{m \in n((\lapha,2(d-1)) - D_n)}{\exp{\{n\phi_{\lapha_n,m/n}\}}} \leq (2(d-1) - \lapha) n^{-K} , 
\end{displaymath}
where $K$ may be chosen to be arbitrarily large by an appropriate choice of $G$. It is this consideration that determines the choice of $G$.

Lemma \ref{lemmanex} implies that the $m$-indexed summand in $C_3(n)$
is at most $r^n \exp{n \phi_{\lapha_n,m/n}}$: thus $C_3(n) \leq (2d+1) r^n$.
Note that $C_1$ satisfies
\begin{eqnarray}
C_1(n) & = & \sum_{m \in n D_n^{*}}{{\bigg( f_n(m/n)\frac{\lapha_n^{m/n}}{(1 + \lapha_n)^{1 + m/n}} \bigg)}^n} \nonumber \\
& \leq &  L n \sum_{m \in n D_n}{g(m/n)^n \exp (n \phi_{\lapha_n,m/n})}, \nonumber 
\end{eqnarray}
where the inequality is a consequence of Theorem \ref{thmone}. The fact that the function $\phi$ is nowhere positive implies that
\begin{displaymath}
 C_1(n) \leq L n \sum_{m \in n D_n}{g(m/n)^n}.
\end{displaymath}
Hence the desired contradiction will be reached if we can show that
\begin{equation}\label{expf}
\sum_{n}{ n \sum_{m \in n D_n}{g(m/n)^n}} 
\end{equation}
is finite.
As such, the proof is completed by the following lemma.
\begin{lemma}
Assume hypotheses $(\reta)$ and $(\appak)$. Suppose that $\appak < 1/2$ and that $\reta \appak < 1$. Then, for $\epsilon \in (0,1 - \reta \appak)$ and $n \in \mathbb{N}$ sufficiently large,
\begin{equation}\label{wrte}
\sum_{m \in n D_n}{g(m/n)^n} \leq \exp{-n^{1 - \reta \appak - \epsilon}}.
\end{equation}
\end{lemma}
{\bf Proof}
Let ${\reta}^* > \reta$ and ${\appak}^* > \appak$ be such that
$\appak^* < 1/2$ and ${\reta}^* {\appak}^* < \reta \appak + \epsilon$.
By hypothesis $(\reta)$, there exists ${\epsilon}'>0$ such that
\begin{displaymath}  
\delta \in (0,{\epsilon}') \ \textrm{implies} \ g(\lapha + \delta) - g(\lapha) < -  {\delta}^{{\reta}^*}.
\end{displaymath}
From Theorems \ref{thmone} and \ref{thmtwo}, it follows that $\sup_{\beta \in [\lapha +
\epsilon',2(d-1)]}{g(\beta)} < 1$, which shows that the contribution
to the sum in (\ref{wrte}) from all those terms indexed by $m$ for
which $m/n > \lapha + \epsilon'$ is exponentially decaying in
$n$. Thus, we may assume that
there exists $N_1$ such that for $n \geq N_1$, if  $m \in D_n^{*}$,
then $m/n - \lapha < \epsilon'$. Note that, by hypothesis $(\appak)$, $\lapha_n - \lapha \geq n^{-{\appak}^*}$ for sufficiently large. Hence, there exists $N_2$ such that, for $n \geq N_2$,
\begin{displaymath}
\lapha_n - G(\log(n)/n)^{1/2} \geq \lapha + n^{-{\appak}^*} - G(\log(n)/n)^{1/2} \geq \lapha + (1/2)n^{-{\appak}^*}.
\end{displaymath}
For $n \geq \max\{ N_1,N_2 \}$ and $m \in n D_n^*$,
\begin{eqnarray}
 g(m/n) & \leq & 1 - (m/n - \lapha)^{{\reta}^*} \nonumber \\
        & \leq & 1 - (\lapha_n - G(\log(n)/n)^{1/2} - \lapha)^{{\reta}^*} \nonumber \\
        & \leq & 1 - ((1/2)n^{-{\appak}^*})^{{\reta}^*}. \nonumber
\end{eqnarray}
So, for $n \geq \max\{N_1,N_2\}$,
\begin{displaymath}
\sum_{m \in n D_n^{*}}{g(m/n)^n} \leq (2G(n \log(n))^{1/2})[1 - C' n^{-{\appak}^* {\reta}^*}]^{n},
\end{displaymath}
for some constant $C' > 0$.
There exists $g \in (0,1)$, such that for large $n$, 
\begin{displaymath}
 [1 - C' n^{-{\appak}^* {\reta}^*}]^{n} \leq g^{n^{1 - {\appak}^*{\reta}^*}}.
\end{displaymath}
This implies that
\begin{displaymath}
\sum_{m \in n D_n^{*}}{g(m/n)^n} \leq  h^{n^{1 - {\appak}^*{\reta}^*}} \ \textrm{for large $n$ and $h \in (g,1)$}.
\end{displaymath}
From ${\reta}^* {\appak}^* < \reta \appak + \epsilon$, we find that
\begin{displaymath}
\sum_{m \in n D_n^{*}}{g(m/n)^n} \leq \exp{-n^{1 - \reta \appak - \epsilon}} \ \textrm{for large $n$,}
\end{displaymath}
as required. $\Box$
\section{Sufficient conditions for $\theta(p_c)=0$}\label{sectthree}
In this section, we prove two theorems, demonstrating sufficient
conditions for the continuity of the percolation probability in terms
of inequalities on $\reta$ and $\appak$.
\subsection{Theorem \ref{thmfour}: when  $\reta < 2$ or  $\appak > 1/2$.}
\begin{theorem}\label{thmfour}
Assume that hypotheses ($\reta$) and ($\appak$) hold.
\flushleft
\begin{enumerate}
\item Suppose that $\reta < 2$. Then $\theta(p_c) = 0$.
\item Suppose that $\appak > 1/2$. Then $\theta(p_c) = 0$.
\end{enumerate}
\end{theorem}
The proof of Theorem \ref{thmfour} will exploit the characterisation of continuity provided by the following lemma.
\begin{definition}
${}$
\flushleft
\begin{itemize}
\item Let $\sigma(p)= \suf{m}{\sigm}$.
\item Let $\sigma_N(p)=\sum_{n \leq N}{\sum_{m}{\sigm}}$
\end{itemize}
\end{definition}
\begin{lemma}\label{lemghe}
A necessary and sufficient condition for $\theta(p_c)=0$ is that $\sigma_n$ tends uniformly to $\sigma$ on the interval $(0,p_c)$.
\end{lemma}
{\bf Proof} 
Note that $\sigma(p)$ is equal to the probability that the origin
belongs to a finite open cluster, and, as such, is identically to one
on $(0,p_c)$. (Note that this requires that we include the possibility that
$C(0)$ is the empty set, which entails setting $\sigma_{0,2d}=1$, and
$\sigma_{0,m}$ for other values of $m$). We have that
\begin{displaymath}
\theta(p_c) = 1 - \sigma(p_c),
\end{displaymath} 
since the sum of critical probabilities that the origin lies in an infinite cluster, or in some finite cluster, is equal to one. This means that the condition $\theta(p_c)=0$ is equivalent to the continuity of $\sigma$ on the closed interval $[0,p_c]$.
The functions $\sigma_n$ are polynomials. Hence, if their uniform limit in $n$ exists on an interval, the limit function is continuous there. The uniform convergence of $\sigma_n$ to $\sigma$ on $(0,p_c)$ and the pointwise convergence of $\sigma_n(p_c)$ to $\sigma(p_c)$ imply that this convergence is uniform on $[0,p_c]$; hence $\sigma_n \to \sigma$ uniformly on $(0,p_c)$ implies $\theta(p_c)=0$.
For the converse, if $\theta(p_c)=0$, then $\sigma(p_c)=1$, and $1 - \sigma_n$ is a decreasing sequence of continuous functions tending to zero pointwise on $[0,p_c]$. Dini's theorem implies that the convergence is uniform. $\Box$ \\  
{\bf Proof of Theorem \ref{thmfour}}
By Lemma \ref{lemghe}, to establish that $\theta(p_c)=0$, it suffices
to show that $\sigma_n$ tends to $\sigma$ uniformly on $(0,p_c)$. We
begin by verifying this condition under the hypotheses of the first
part of the Theorem. We will show that 
\begin{equation}\label{escc}
\sum_{n}{\sum_{m}{\sigma_{n,m}\sup_{p \in (0,p_c)}p^n(1-p)^m}} < \infty .
\end{equation}
This will do because
\begin{eqnarray}
\sup_{p \in (0,p_c)}{\big( \sigma(p) - \sigma_{N}(p) \big)} & = & \sup_{p \in (0,p_c)}{\sum_{n \geq N+1}{\sum_{m}\sigm}} {}\nonumber\\
 & \leq & \sum_{n \geq N+1}{\sum_{m}{\sdash \sup_{p \in (0,p_c)}{p^n(1-p)^m}}} {}\nonumber
\end{eqnarray}
So the condition (\ref{escc}) implies the uniform convergence of $\sigma_n$ to $\sigma$ on the subcritical interval.

Note that
\begin{displaymath}
\sup_{p \in (0,p_c)}{p^n(1-p)^m} = \left\{ \begin{array}{ll} \sigt & \textrm{for $n/(n+m) \leq p_c$} \\
\sigo & \textrm{for other pairs $(n,m)$}. \end{array} \right..
\end{displaymath}

This observation allows us to decompose the sum appearing in (\ref{escc}):
\begin{eqnarray}
\suf{m}{\sdash \sup_{p \in (0,p_c)}{p^n(1-p)^m}} & = & \suh{m  = 1}{\lfloor n \lapha \rfloor }{\sigma_{n,m} \sigo } {}\nonumber\\
& + & \suf{m > \lfloor n \lapha \rfloor}{\sigma_{n,m}\sigt} .{}\nonumber
\end{eqnarray}
Now,
\begin{displaymath}
 \suh{m =  1}{\lfloor n \lapha \rfloor}{\sdash \sigo } \leq \suf{m}{\sdash \sigo},
\end{displaymath}
which is less than or equal to one, being the critical probability that the origin lies in a finite cluster.

Set
\begin{displaymath}
A = \suf{m > \lfloor n \lapha \rfloor}{\sdash\sigt}.
\end{displaymath}
Then it suffices to show that $A$ is finite. Our strategy is to split each of the summands of $n$ into two parts, each of which is a sum over $m$ in an interval which has an $n$-dependence. The first sum, $A_1$, will include those $m$-values sufficiently close to $n \lapha$ that this term can be bounded in terms of the critical probability of observing a large cluster. The second sum, $A_2$, will be shown to decay quickly, under the assumption that $\reta < 2$.

Write $A = A_1 + A_2$, where
\begin{eqnarray} 
A_1 & = & \suh{m =  \lfloor n \lapha \rfloor + 1}{\lfloor n \lapha + n^{1/2} \rfloor +1 }{\sdash\sigt}, \nonumber \\
\textrm{and} \, \, A_2 & = & \suf{m > \lfloor n \lapha + n^{1/2} \rfloor +1}{\sdash\sigt}.\nonumber
\end{eqnarray}
Recalling that $\lapha = 1/p_c - 1$,
\begin{displaymath}
A_1 =  \suh{m = \lfloor n \lapha \rfloor + 1}{\lfloor n
\lapha + n^{1/2} \rfloor + 1 }{\sdash\sigo \exp( - n \phi (\lapha,m/n)
)},  
\end{displaymath}
where the function $\phi$ was given in Definition \ref{defnphi}.
For each $m \in \{ \lfloor n \lapha \rfloor, \ldots, \lfloor n \lapha
+ n^{1/2} \rfloor + 1 \}$, $c_m \in (0,3/2)$, where $c_m$ is given by $m/n = \lapha + c_m n^{-1/2}$. The quantities $c_m$ could be chosen to lie in any open interval of the form $(0, 1+ \epsilon)$; a choice has been made for definiteness.

Lemma \ref{lemfour} implies that for any sufficiently large $C'$, there exists $N_1$ such that for all $n \geq N_1$, and for $m \in \{ \lfloor n \lapha \rfloor + 1, \ldots, \lfloor n \lapha + n^{1/2} \rfloor + 1 \}$,
\begin{displaymath}
 - \phi (\lapha,m/n) \leq 9/[8 n \lapha(\lapha +1)] + C'/{n^{3/2}} .
\end{displaymath}
From this, we deduce that for $n \geq N_1$ and $m \in \{
\lfloor n \lapha \rfloor + 1, \ldots, \lfloor n \lapha + n^{1/2}
\rfloor + 1 \}$, \mbox{$\exp{(- n \phi (\lapha,m/n) )}$} is bounded
above, by $C$, say.
So,
\begin{eqnarray}
A_1 & \leq & \sum_{n < N_1}{\sum_{m \in \{ \lfloor n \lapha \rfloor + 1, \ldots, \lfloor n \lapha + n^{1/2} \rfloor + 1 \}}{\sdash\sigt}} {} \nonumber \\
   {} & & + C   \sum_{n \geq N_1}{\sum_{m \in \{ \lfloor n \lapha \rfloor + 1, \ldots, \lfloor n \lapha + n^{1/2} \rfloor + 1 \}}{\sdash\sigo}}, \nonumber 
\end{eqnarray}
which is finite, as desired.

We now seek to bound $A_2$:
\begin{eqnarray}
A_2 & = & \suf{\atwo}{{\bigg(f_n(m/n)\Big(\frac{n}{n+m}\Big){\Big(\frac{m}{n+m}\Big)}^{m/n}\bigg)}^n} {} \nonumber \\          
    & \leq & L \sum_{n}{\sum_{m = \lfloor n \lapha + n^{1/2} \rfloor +
    2}^{2(d-1)n - 1}{  n
    \bigg(f(m/n) \Big(\frac{n}{n+m}\Big){\Big(\frac{m}{n+m}\Big)}^{m/n}\bigg)^n }} \nonumber \\
   & & + {}  \sum_{n}{\sum_{m = 2(d-1)n}^{2(d-1)n + 2d}{\sdash\sigt}} \nonumber
\end{eqnarray}
where the inequality follows from Theorem \ref{thmone} and the fact
that $g \leq 1$. By Lemma \ref{lemmanex}, there exists $r \in (0,1)$ such that, for $n$ sufficiently large,
\begin{displaymath}
 \sum_{m = 2(d-1)n}^{2(d-1)n + 2d}{\sdash\sigt} \leq (2d+1) r^n . 
\end{displaymath}
It follows from the definition of the function $g$ that
\begin{equation}\label{eat}
A_2 \leq L  \sum_{n}{\sum_{m = \lfloor n \lapha + n^{1/2} \rfloor + 2}^{2(d-1)n - 1}{n g(m/n)^n}} + (2d+1) \sum_{n}{r^n}.
\end{equation}
We aim to bound $g$ above on intervals of the form $(\lapha,\lapha +
\epsilon)$, and to use these bounds to show that the first term in the expression on the
right-hand-side of [\ref{eat}] is finite. The condition $\reta < 2$ enables us to do this.
Let $\epsilon \in (0,2 - \reta)$. Let ${\delta}' > 0$, be such that, for $\delta \in (0,{\delta}')$, $g(\lapha + \delta) - g(\lapha) < -{\delta}^{\reta + \epsilon}$. Let $\gamma \in (0,1)$ be such that 
\begin{displaymath}
\sup_{\beta \in (\lapha + {\delta }',2(d-1))}{g(\beta)} < \gamma .
\end{displaymath}
Note that 
\begin{eqnarray}
     &  & \suh{m = \lfloor n \lapha + n^{1/2} \rfloor
     +2}{\lfloor n(\lapha + {\delta}') \rfloor}{ n g(m/n)^n} {} \nonumber \\
 & \leq & \suh{m = \lfloor n \lapha + n^{1/2} \rfloor
 +2}{\lfloor n(\lapha + {\delta}') \rfloor}{ n \Big( 1 - \big( m/n -
 \lapha \big)^{\reta + \epsilon} \Big)^n} {} \nonumber \\
    & \leq & {\delta}' \sum_{n}{n^2 \Big( 1 - n^{ - \frac{\reta +
    \epsilon}{2}} \Big)^n}. \nonumber 
\end{eqnarray}
Since $\reta + \epsilon < 2$, this expression is finite.
Note also that
\begin{eqnarray}
 & &  \suh{m  = \lfloor n(\lapha + {\delta}') \rfloor +
 1}{2(d-1)n - 1}{n g(m/n)^n}  \nonumber \\
 & \leq & 2(d-1)\sum_{n}{n^2 \gamma^n} < \infty . 
\end{eqnarray}
We deduce that $A_2$ is finite and in doing so, complete the proof of
the first part of Theorem \ref{thmfour}.

We now prove the second part of the Theorem.
A sufficient condition for continuity is
\begin{equation}\label{anothscc}
\suf{m}{\sdash t_n^n (1 - t_n)^m} < \infty .
\end{equation}
Indeed, the supremum over $p$ in $(0,p_c)$ of $\sigma - \sigma_N$ is
bounded above by the expression in (\ref{anothscc}) with the sum in $n$
being taken over values exceeding $N - 1$. By Lemma \ref{lemghe}, if
(\ref{anothscc}) holds, then $\theta(p_c)=0$.

We will decompose the sum that appears in (\ref{anothscc}) and try
to bound the various parts. Several of these parts may be shown to be
finite without invoking the two hypotheses, $(\reta)$ and
$(\appak)$. Firstly, we present the arguments involving these parts;
then, we show how the condition $\appak > 1/2$ implies finiteness of
the sum in (\ref{anothscc}).

We have that
\begin{eqnarray}
\suf{m}{\sdash t_n^n(1 - t_n)^m} & = & \suh{m = 1}{\lfloor n \lapha \rfloor}{\sdash t_n^n(1 - t_n)^m} \nonumber \\
  & + & \suh{m = \lfloor n \lapha \rfloor + 1}{2(d-1)n + 2d}{\sdash t_n^n(1 - t_n)^m}  \nonumber 
\end{eqnarray}
For given $n,m \in \mathbb{N}$, the function on $[0,1]: t \to t^n(1-t)^m$ attains its supremum at $n/(n+m)$, and is increasing on $[0,n/(n+m)]$, and decreasing on $[n/(n+m),1]$.

The fact that $t_n \leq p_c$ implies that $t_n^n(1 - t_n)^m \leq p_c^n(1 - p_c)^m$ provided that $n/(n+m) > p_c$, which holds if and only if $m \leq \lfloor n \lapha \rfloor$. This observation enables us to bound the first of the two sums:
\begin{eqnarray}\label{wret} 
 \suh{m = 1}{\lfloor n \lapha \rfloor}{\sdash t_n^n(1 - t_n)^m} 
 & \leq & \suh{m = 1}{\lfloor n \lapha \rfloor}{\sdash \sigo} \nonumber \\
& \leq &  \suf{m}{\sdash \sigo} \leq 1 \nonumber
\end{eqnarray}
It remains to analyse the expression 
\begin{equation}
\suh{m = \lfloor n \lapha \rfloor + 1}{2(d-1)n + 2d}{\sdash t_n^n(1 - t_n)^m}. 
\end{equation}
To do so, we make the following  definition.
\begin{definition}
For $G \in \mathbb{N}$, let  $D_n^{*} ( = D_n^{*}(G))$ denote the interval 
\begin{displaymath}
D_n^{*} = (\max{\{ \lapha, \lapha_n - G{\{\log(n)/n \}}^{1/2} \}}, \lapha_n + G{\{\log(n)/n \}}^{1/2}),
\end{displaymath} 
where the constants $\{ \lapha_n : n \in \mathbb{N} \}$ were specified
in Definition \ref{alphn}.
\end{definition}
Allowing that $G$ will be determined slightly later, we write the
expression in (\ref{wret}) in the form
\begin{eqnarray}
&  & \suf{m \in n D_n^{*}}{\sdash t_n^n(1 - t_n)^m} + \suf{m \in n((\lapha,2(d-1)) - D_n^{*})}{\sdash t_n^n(1 - t_n)^m} \nonumber \\
& & \, + \,  \sum_{n}{\sum_{m = 2(d-1)n}^{2(d-1)n + 2d}{\sdash t_n^n(1
- t_n)^m}} . \label{fvg}  \nonumber
\end{eqnarray}
An argument identical to that by which the term $C_2$ was bounded in the
proof of Theorem \ref{thmthr} yields
\begin{displaymath}
 \suf{m \in n((\lapha,2(d-1)) - D_n^{*})}{\sdash t_n^n(1 - t_n)^m}
 \leq \sum_{n}{n^{-K}}, 
\end{displaymath}
where $K$ may be chosen to be arbitrarily large by an appropriate
choice of $G$, thereby determining how $G$ is chosen. The third term
in (\ref{fvg}) was labelled $C_3(n)$ in the proof of Theorem \ref{thmthr}
and was shown to be bounded above by $(2d+1) r^n$ for $n$ sufficiently
high. As for the first, we have that
\begin{eqnarray}
& & \suf{m \in n D_n^{*}}{\sdash t_n^n(1 - t_n)^m} \nonumber \\
& = & \suf{m \in n D_n^{*}}{\sdash  \frac{\lapha_n^m}{(1 + \lapha_n)^{n+m}
}} \nonumber \\
& = & \suf{m \in n D_n^{*}}{\sdash \frac{\lapha^m}{(1 + \lapha)^{n+m}}
 \exp{n\Phi(\lapha_n,\lapha,m/n)} }, \label{kzc}
\end{eqnarray}
where 
\begin{eqnarray}
\Phi(\gamma,\lapha,\beta) & = & \beta \log \gamma - (\beta + 1) \log ( \gamma + 1) - \beta \log \lapha + (\beta + 1)\log(\lapha + 1) \nonumber \\
  & = & \beta \log (1 + (\gamma - \lapha)/\lapha ) - (\beta + 1) \log (1 + (\gamma - \lapha)/(1+ \lapha)) \nonumber \\
  & = & - \, \,  \frac{(\gamma - \lapha)^2}{2} \bigg[
  \frac{\beta}{{\lapha}^2} -  \frac{\beta + 1}{(1 + \lapha)^2} \bigg] \nonumber \\
 & &  + \, \, \frac{(\gamma - \lapha)(\beta - \lapha)}{\lapha(\lapha + 1)} \,  + \, O\big[(\gamma - \lapha)^3 \big]. \nonumber
\end{eqnarray}
We are supposing that hypothesis $(\appak)$ holds, and that $\appak > 1/2$. Let ${\appak}'$ satisfy $\appak > {\appak}' > 1/2$.
In this context,
\begin{eqnarray}
\Phi(\lapha_n,\lapha,\beta) & = &  \frac{(\lapha_n - \lapha)(\beta -
\lapha)}{\lapha(\lapha + 1)} \nonumber \\
 & & - \  \frac{(\lapha_n - \lapha)^2}{2} \bigg[
 \frac{\beta}{{\lapha}^2} - \frac{\beta + 1}{(1 + \lapha)^2} \bigg] +
 O \big( n^{-3{\appak}'} \big) . \nonumber
\end{eqnarray}
Now, $\beta \in D_n^{*}$ implies that there exists $C' > 0$ such that $\beta - \lapha \leq C'n^{-{\appak}'} +C'(\log(n)/n)^{1/2}$ ; since $\appak' > 1/2$, we may write $\beta - \lapha \leq C'(\log(n)/n)^{1/2}$, where the value of $C'$ has been increased if necessary. For such $\beta$, $\Phi(\lapha_n,\lapha,\beta) \leq C'n^{-{\appak}' - {1/2}}{\log(n)}^{1/2} + n^{-2{\appak}'} + O(n^{-3{\appak}'})$. This implies that, for all $n$ and $\beta \in D_n^{*}$, $\exp{n\Phi(\lapha_n,\lapha,\beta)} < C'$, where once again the value of $C'$ may have changed.
Recalling that $\lapha = 1/p_c - 1$, we deduce from (\ref{kzc}) that
\begin{eqnarray}
 & & \suf{m \in n D_n^{*}}{\sdash t_n^n(1 - t_n)^m} \nonumber \\
 & \leq & C' \suf{m \in n D_n^{*}}{\sdash \sigo} \nonumber \\
 & \leq & C' \suf{m}{\sdash \sigo} \leq C', \nonumber
\end{eqnarray}
proving the second part of Theorem \ref{thmfour}.$\Box$

\subsection{Theorem \ref{thmund}: when $\appak < 1/2$ and $\reta \appak > 1$.}
We now examine the case where
 $\appak < 1/2$ and $\reta \appak > 1$.
\begin{definition}
Let $n \in \mathbb{N}$, and $\beta \in (0,2(d-1))$. Set 
\begin{equation}\label{eanb}
a_n(\beta) = {\bigg( \frac{f_n(\beta)}{f(\beta)} \bigg) }^n.
\end{equation}
\end{definition}
{\bf Remark} The quantities $a_n(\beta)$ appear in the factorisation of $\sigma_{n,\lfloor \beta n \rfloor}$,
\begin{equation}\label{zcx}
 \sigma_{n,\lfloor \beta n \rfloor} = a_n(\beta)g(\beta)^n \bigg( \frac{(\beta + 1)^{\beta + 1}}{\beta^\beta} \bigg)^n;
\end{equation}
the number of animals containing the origin with given surface-area-to-volume ratio grows at an exponential rate which is conveniently represented as a fraction of an explicit form. The extent to which that rate is underestimated is expressed by the $a_n(\beta)$.
 
We will perform a similar analysis to those undertaken in Theorem \ref{thmfour} to obtain the following result. 
\begin{theorem}\label{thmund}
Assume that hypotheses $(\reta)$ and $(\appak)$ hold. Suppose that \mbox{$\appak < 1/2$} and \mbox{$\reta \appak > 1$}. Let $K$ be large. Then there exist constants $\epsilon > 0$ and $C > 0$ such that for each $n \in \mathbb{N}$,
\begin{eqnarray}
 \epsilon \sum_{\alplow}{a_n(m/n)} & \leq & \sum_{m}{\sdash \sigo} \label{edf} \\
                                   & \leq & \sum_{\alphigh}{a_n(m/n)}
				   + \ n^{-K} \nonumber
\end{eqnarray}
and
\begin{eqnarray}
 \epsilon \sum_{\alpnlow}{a_n(m/n)} & \leq & \sum_{m}{\signt} \label{edg} \\
                                    & \leq  & \sum_{\alpnhigh}{a_n(m/n)} + \ n^{-K} . \nonumber 
\end{eqnarray}
\end{theorem}
{\bf Remark}  Here, $B(a,b)$ denotes the interval $(a-b,a+b)$. \\
Before embarking on the proof of the theorem, we state and prove two
lemmas. 
\begin{lemma}\label{lemest}
For any $K > 0$, there exist constants $\epsilon > 0$ and $C > 0$ such that
\begin{itemize}
\item $\exp{\phi(\lapha,\beta)} > \epsilon$ for $\beta \in (\lapha - n^{-1/2},\lapha + n^{-1/2})$ 
\item  $\exp{\phi(\lapha_n,\beta)} > \epsilon$ for $\beta \in (\lapha_n - n^{-1/2},\lapha_n + n^{-1/2})$ 
\item  $\exp{\phi(\lapha,\beta)} < n^{-K}$ for $\beta \in (0,2(d-1)) -
B \big( \lapha, C (\log n / n)^{1/2} \big) $
\item  $\exp{\phi(\lapha_n,\beta)} < n^{-K}$ for $\beta \in (0,2(d-1)) - 
B \big( \lapha_n, C (\log n / n)^{1/2} \big) $ 
\end{itemize}
\end{lemma}
{\bf Proof}
Each statement follows from the fact there exists $H$ such that $\phi$ satisfies 
\begin{displaymath}
\phi(\lapha,\beta) = -H(\beta - \lapha)^2 + O((\beta - \lapha)^3) \ \textrm{as $\beta \to \lapha$}. \ \Box
\end{displaymath}

In the second of the two lemmas, we invoke the inequalities on $\reta$ and $\appak$ in order to provide a bound on $g(\beta)$ on the range of values of $\beta$ we are considering.

\begin{lemma}\label{lemglb}
Assume that hypotheses $(\reta)$ and $(\appak)$ hold. Suppose that \mbox{$\appak < 1/2$} and \mbox{$\reta \appak > 1$}. Then, for $n$ sufficiently large,
\begin{displaymath}
g(\lapha_n + \delta)^n \geq \frac{1}{2} \ \textrm{for} \ \delta \in (-n^{-1/2},n^{-1/2}).
\end{displaymath}
\end{lemma}
{\bf Proof} 
Let ${\appak}^*,{\reta}^*$ satisfy $\appak < {\appak}^* < 1/2$, $2 < {\reta}^* < \reta$, and ${\reta}^* {\appak}^* > 1$.
Hypothesis $(\reta)$ implies that for any ${\reta}^* \in (2,\reta)$, there exists ${\delta}' > 0$ such that, for $\delta \in (0,{\delta}')$, $g(\lapha + \delta) > 1 - {\delta}^{{\reta}^*}$.
For $\delta \in (-n^{-1/2},n^{-1/2})$ and $n$ sufficiently large, we have that
\begin{eqnarray}
g(\lapha_n + \delta) & \geq & 1 - (\lapha_n + \delta - \lapha)^{{\reta}^*} \ \textrm{and} \nonumber \\
\lapha_n + \delta - \lapha & \leq & n^{-{\appak}^*} + n^{-1/2} \leq 2n^{-{\appak}^*} \nonumber
\end{eqnarray}
For such $n$,
\begin{displaymath}
g(\lapha_n + \delta) \geq 1 - 2^{{\reta}^*}n^{-{\appak}^*{\reta}^*},
\end{displaymath}
which implies that
\begin{displaymath}
g(\lapha_n + \delta)^n \geq 1/2 \ \textrm{for} \ \delta \in (-n^{-1/2},n^{-1/2}) \ \textrm{and $n$ sufficiently large.} \, \Box
\end{displaymath}
{\bf Proof of Theorem \ref{thmund}.}
An important element in the proof is expressed in the representation:
\begin{eqnarray}
 & & \sum_{m = 1 }^{2(d-1)n - 1}{\sdash \sigo} \nonumber \\
 & = &  \sum_{m = 1}^{2(d-1)n - 1}{{\bigg( f_n(m/n)\frac{\lapha^{m/n}}{(1 + \lapha)^{1 + m/n}} \bigg)}^n} \nonumber \\
 & = & \sum_{m = 1}^{2(d-1)n - 1}{a_n(m/n)g(m/n)^n \exp{n \phi_{\lapha,m/n}}}, \nonumber 
\end{eqnarray}
where the function $\phi$ was specified in Definition \ref{defnphi}.
Similarly,
\begin{displaymath}
\sum_{m}{\signt} =  \sum_{m}{a_n(m/n)g(m/n)^n \exp{n \phi_{\lapha_n,m/n}}}.
\end{displaymath}
Firstly, we compute the lower bounds in the  sums (\ref{edf}) and (\ref{edg}). By Lemmas \ref{lemest} and \ref{lemglb}, there exists $\epsilon > 0$, such that,  for $n$ sufficiently large,
\begin{eqnarray}
 & & \qquad \quad \sum_{m}{\sigma_{n,m}p_c^n(1 - p_c)^m} \nonumber \\
 & \geq & \sum_{m = \lfloor n \lapha \rfloor +1}^{\lfloor  n(\lapha + n^{-1/2}) \rfloor}{a_n(m/n)g(m/n)^n \exp{n \phi_{\lapha,m/n}}} \nonumber \\
 & + &  \sum_{m = \lfloor n(\lapha - n^{-1/2}) \rfloor + 1}^{\lfloor n \lapha \rfloor}{a_n(m/n)g(m/n)^n \exp{n \phi_{\lapha,m/n}}} \nonumber \\
 & \geq & {\epsilon}^2 \sum_{m = \lfloor n \lapha \rfloor + 1}^{\lfloor n(\lapha + n^{-1/2}) \rfloor}{a_n(m/n)} + \epsilon \sum_{m = \lfloor n(\lapha - n^{-1/2}) \rfloor + 1}^{\lfloor n \lapha \rfloor}{a_n(m/n)}.\nonumber 
\end{eqnarray}
Similarly,
\begin{eqnarray}
 \sum_{m}{\signt} & \geq & \sum_{\alpnlow}{a_n(m/n)g(m/n)^n \exp{n \phi_{\lapha,m/n}}} \nonumber \\
 & \geq & \epsilon \sum_{\alpnlow}{a_n(m/n)g(m/n)^n} \nonumber \\
 & \geq & {\epsilon}^2 \sum_{\alpnlow}{a_n(m/n)}. \nonumber
\end{eqnarray}
To estimate the upper bounds, we proceed as follows.
\begin{displaymath}
\sum_{m}{\sdash \sigo} = D_1 + D_2 + D_3,
\end{displaymath}
where
\begin{eqnarray}
D_1 & = & \sum_{\alphigh}{a_n(m/n)g(m/n)^n  \exp{n \phi_{\lapha,m/n}}}, \nonumber \\
D_2 & = & \sum_{m \in n E_n }{a_n(m/n)g(m/n)^n  \exp{n
\phi_{\lapha,m/n}} }, \nonumber 
\end{eqnarray}
and
\begin{displaymath}
 D_3  =  \sum_{m = 2(d-1)n }^{2(d-1)n + 2d }{\sdash \sigo} . 
\end{displaymath}
Here, the constant $C$ is determined by Lemma \ref{lemest}, while
$E_n$ denotes the set $(0,2(d-1)) - B(\lapha, C(\log n / n )^{1/2} )$.
We have that
\begin{displaymath}
D_1 \leq  \sum_{\alphigh}{a_n(m/n)},
\end{displaymath}
and, by Lemma \ref{lemest},
\begin{eqnarray}
D_2 & \leq & L n \sum_{m \in n E_n}{\exp{n \phi_{\lapha,m/n}}} \nonumber \\ 
    & \leq & 2(d-1)L n^{1-K}. \nonumber 
\end{eqnarray}
Note also that by Lemma \ref{lemmanex}, for some $r \in (0,1)$ and $n$ sufficiently large,
\begin{displaymath}
D_3 \leq (2d + 1) r^n  .
\end{displaymath}
Similarly, for $n \in \mathbb{N}$,
\begin{displaymath}
\sum_{m}{\signt} = D_1^{*} + D_2^{*} + D_3^{*},
\end{displaymath}
where
\begin{eqnarray}
D_1^{*} & = & \sum_{\alpnhigh}{a_n(m/n)g(m/n)^n \exp{n \phi_{\lapha_n,m/n}}}, \ \nonumber \\
D_2^{*} & = & \sum_{m \in n E'_n }{a_n(m/n)g(m/n)^n \exp{n
\phi_{\lapha_n,m/n}}}, \nonumber 
\end{eqnarray}
and
\begin{displaymath}
 D_3^{*}  =  \sum_{m = 2(d-1)n}^{2(d-1)n + 2d}{\signt},
\end{displaymath}
where 
$E'_n$ denotes the set $(0,2(d-1)) - B(\lapha_n, C(\log n / n )^{1/2} )$.
For each $n \in \mathbb{N}$,
\begin{displaymath}
D_1^{*} \leq  \sum_{\alpnhigh}{a_n(m/n)},
\end{displaymath}
and, by Lemma \ref{lemest},
\begin{eqnarray}
D_2^{*} & \leq & L n \sum_{m \in n E'_n}{\exp{n \phi_{\lapha_n,m/n}}} \nonumber \\ 
    & \leq & 2(d-1)Ln^{1-K}. \nonumber 
\end{eqnarray}
Finally, by Lemma \ref{lemmanex},
\begin{displaymath}
D_3^{*} \leq (2d+1) r^n .
\end{displaymath}  
This completes the proof of the theorem. $\Box$ \\
{\bf Remark} From the proof of Theorem \ref{thmund}, we see that the
condition 
\begin{displaymath}
\sum_{\alpnhigh}{a_n(m/n)} < \infty
\end{displaymath}
implies that $\theta(p_c)=0$, without recourse to scaling
hypotheses. In examining this condition, bounds on the entropic
exponent are revelant (see \cite{MR95m:82076}).  

\subsection{Related exponent inequalities}

In \cite{MR88k:60176}, C.M. Newman proved a sufficient condition for the absence of an infinite cluster at criticality.
We reproduce the theorem of that paper, and its brief proof, making some changes to reconcile the notation with that of this report.

\begin{theorem}{[Newman]}
For Bernoulli bond percolation on $\mathbb{Z}^d$, define the mean
cluster size $\chi(p)$ for $p \in (0,p_c)$ as
\begin{displaymath}
\chi(p)= \sum_{n}{n \mathbb{P}_p(\vert C(0) \vert = n)}.
\end{displaymath}   
If $\int_{0}^{p_c}{\chi(p)^{1/2} dp} < \infty$, then $\theta(p_c)=0$.
\end{theorem}
{\bf Remark} The exponent for mean cluster size $\gamma$ is given by $\lim_{p \uparrow
p_c}{\log \chi(p) / \log p}$, if this limit exists. Newman's result
could be phrased: if $\gamma < 2$, then $\theta(p_c)=0$. \\ 
{\bf Proof.}
The probability $\mathbb{P}_p(\vert C(0) \vert = n )$ is given by
\begin{equation}\label{eqa}
\mathbb{P}_p(\vert C(0) \vert = n ) = \sum_{m}{\sigma_{n,m}p^n(1-p)^m},
\end{equation}
and the percolation probability $\theta$ may be written,
\begin{equation}\label{eqb}
\theta(p) = 1 - \sum_{n}{\mathbb{P}_p(\vert C(0) \vert = n )}.
\end{equation}
Substituting  (\ref{eqa}) in (\ref{eqb}) and  differentiating gives that
\begin{equation}\label{eqc}
\sum_{n}{\sum_{m}{\big( n/p - m/(1-p) \big) \sigma_{n,m} p^n (1-p)^m }} = 0.
\end{equation}
Differentiating again gives that
\begin{eqnarray}
& & \sum_{n}{\sum_{m}{\big( n/p - m/(1-p) \big)^2 \sigma_{n,m} p^n (1-p)^m }} \nonumber \\
& = & \sum_{n}{\sum_{m}{\big( n/{p^2} + m/{(1-p)^2} \big) \sigma_{n,m} p^n {(1-p)}^m}} \nonumber \\
& = & \frac{\chi(p)}{p^2 (1-p)} \ \textrm{for $0 < p < p_c$,}
\end{eqnarray}
where (\ref{eqc}) was used in the last inequality.
The term-by-term differentiation was justified by the exponential decay rate of $\mathbb{P}_p(\vert C(0) \vert = n )$ for $p < p_c$. (cf \cite{MR86h:82045}).
For $\epsilon > 0$ small, we have that
\begin{eqnarray}
\theta(p_c) & = & \theta(p_c) - \theta(p_c - \epsilon) \nonumber \\
            & = & \lim_{n \to \infty}{\int_{p_c - \epsilon}^{p_c}{\sum_{n < N}{-\frac{d}{dp}\mathbb{P}_p(\vert C(0) \vert = n )dp}}} \nonumber \\
    & = & \lim_{n \to \infty}{\int_{p_c - \epsilon}^{p_c}{\sum_{n < N}{\sum_{m}{- \big ( n/p - m/(1-p) \big) \sigma_{n,m} p^n (1-p)^m }}}} \nonumber \\
   & \leq &  \lim_{n \to \infty}{\int_{p_c - \epsilon}^{p_c}{ \bigg( \sum_{n < N}{\sum_{m}{ \big ( n/p - m/(1-p) \big)^2 \sigma_{n,m} p^n (1-p)^m }} \bigg)^{1/2} }} \nonumber \\
 & = & \int_{p_c - \epsilon}^{p_c}{\bigg( \frac{\chi(p)}{p^2 (1-p)} \bigg)^{1/2} dp}, \label{eqd}
\end{eqnarray}
where the last inequality is a consequence of the Cauchy-Schwarz inequality on the sequence space $\{f_{n,m}: \sum_{n}{\sum_{m}{{\vert f_{n,m} \vert}^2} \sigma_{n,m} p^n (1-p)^m } < \infty \}$. 
The result follows immediately from (\ref{eqd}). $\Box$\\

The exponent $\delta$ describes the decay rate of the probability of
large clusters at the critical value. If it exists, then
$\mathbb{P}_{p_c}(\vert C(0) \vert = n )$ decays at the rate $n^{-(1 + 1/{\delta})}$.
Theorem $1.3$ of \cite{MR88k:60176} asserts that $\gamma \geq 2(1 -
1/{\delta})$. Its proof has the flavour of the arguments in Theorem
(\ref{thmfour}).

\subsection{Underestimating $f$: rigorous results and conjecture}
Recall the factorisation of $\sigma_{n,m}$ that appears in (\ref{zcx}).
The quantities $a_n(\beta)$ measure the extent of the miscalculation of the asymptotic exponential rate, $f(\beta)$, by the computed exponential rate, $f_n(\beta)$. How do the $a_n(\beta)$ behave for large $n$, in different $\beta$-intervals?
A little light is shed on this question by work of N. Madras \cite{MR2001f:82034}. Setting $a_n$ equal to the number of lattice animals up to translation with $n$ vertices, and writing $\mu = \lim_{n \to \infty}{a_n^{1/n}}$, the author makes the following definition:
\begin{definition}
Let the sequence $\{ \theta_n : n \in \mathbb{N} \}$ be given by the relation
\begin{displaymath}
a_n = n^{- \theta_n} \mu^n.
\end{displaymath}
\end{definition}
His paper demonstrates that $\limsup_{n \to \infty}{\theta_{n}} \geq 1 - 1/d$. The argument involves concatenating any pair of $n$-edged animals to form a new animal, from which the pair may be recovered. Madras points out that the second animal may be joined to the first in a variety of different positions. The lexicographically smallest vertex of the second animal may be translated to any position neighbouring a vertex of the first animal on a `face' of that animal with the most vertices, before the joining takes place. It is this surface-area effect which enables the term $1 - 1/d$ to appear in the inequality.   
Madras' argument may be adapted and unified with the correction construction proof of Theorem \ref{thmone} to yield the following inequality:

\begin{displaymath}
\textrm{for $\beta \in (0,2(d-1))$,} \ \limsup_{n \to \infty}{\frac{- \log a_n(\beta)}{\log n}} \geq -1/d.
\end{displaymath}

There are order-$n$ as many $n$-edged clusters containing the origin
as there are $n$-edged animals up to translation. Thus, the change in
the right-hand-side from $1 - 1/d$.

Other work is relevant for understanding the behaviour of the
quantities $\{ a_n(\beta) : n \in \mathbb{N} \}$. Unusually large finite clusters occur in the supercriticial phase of bond percolation as a result of a surface-area effect, in which a large cluster of edges with the characteristics of a region of the infinite cluster happens to be cut off from that cluster by a collection of closed edges which form a surface that encloses the finite cluster. Cerf (see \cite{MR1774341}) has determined the decay rate of the probability that the cluster to which the origin belongs has more than $n$ edges in the three-dimensional bond percolation model. The following result is a consequence of Theorem $2.5$ of that work.

\begin{theorem}\label{thmcerf}
In the three dimensional bond percolation model for $p > p_c$, there
exists a constant $D (=D(p)) > 0$ such that
\begin{displaymath}
\lim_{n \to \infty}{\frac{\log \mathbb{P}_p(n^3 \leq \vert C(0) \vert < \infty)}{n^2}} = D
\end{displaymath}
\end{theorem}
Cerf determines the constant $D$. It may be expressed in terms of the surface tension of a convex body, called the Wulff crystal, which minimises a variational problem; this body is the likely shape of an unusually large finite cluster in the supercritical regime, after the gaps in its internal structure (present on a `mesoscopic' scale) have been eliminated.

It is possible to demonstrate by an argument similar to that used in the proof of Theorem \ref{thmund}, that the following result holds.

\begin{theorem}
For $p> p_c$, and for any $C' > 0$, we may
find a small positive $\epsilon$ and a large positive $C$ such that,
\begin{eqnarray}
 & & \epsilon \sum_{m \in n B(1/p - 1, c n^{-1/2} ) }{a_n(m/n)}
 \nonumber \\
& \leq & \mathbb{P}_p(\vert C(0) \vert = n) \nonumber \\ 
& \leq & \sum_{m \in n B(1/p -1, C n^{-1/2d}) }{a_n(m/n)} \, + \,  \exp{- C' n^{1 - 1/d}}, \nonumber
\end{eqnarray}
\end{theorem}
Combining this result with Cerf's work suggests the conjecture, \\
\begin{conjecture}
Let $d=3$. For $\beta \in (0,\lapha)$, the limit
\begin{displaymath}
\lim_{n \to \infty}{\frac{- \log a_n(\beta)}{n^{2/3}}}
\end{displaymath}
exists, and is equal to the constant $D(1/(1+ \beta))$ that appears in
Theorem \ref{thmcerf}.
\end{conjecture}

There appears to be no rigorous work on the decay rate of $a_n(\beta)$ for fixed $\beta$-values lying in the interval $(\lapha,2(d-1))$, beyond the adaptation of Madras' concatenation argument. Meir, Aharony and A. Brooks Harris (see \cite{MR89k:82058}) perform numerical studies which are relevant. 
Their work would lead one to suppose that
for $\beta \in (\lapha,2(d-1))$, the limit $\lim_{n \to \infty}{\frac{- log \, a_n(\beta)}{\log n}}$ exists, and is independent of $\beta$ on this interval.
This belief is consistent with the idea of universality. In this case, the hypothesis of universality might be expressed by saying that lattice animals with surface-area-to-volume ratio $\beta$ lie in the same universality class for the entropic exponent, for $\beta \in (\lapha,2(d-1))$.

\section{Scaling law}\label{sectfour}
In this section, we examine the exponential decay rate in $n$ for the
probability of the event $\{ C(0)=n \}$ for $p$ slightly less than
$p_c$ by our combinatorial approach. In doing so, we relate the quantity $\reta$ to the exponent
for correlation size, and see how the scaling behaviour for the typical
surface-area-to-volume ratio of unusually large clusters in the
marginally subcritical regime depends on the value of $\reta$.   
\subsection{Decay rate for large subcritical clusters}
\begin{definition}
Let $\xi: (0,1) \times (0,\infty) \to \mathbb{R}$ be given by
\begin{displaymath}
\xi(p,\beta) = \beta \log \beta - (\beta + 1) \log(\beta + 1) -
\log{p} - \beta \log(1-p).
\end{displaymath}
\end{definition}
Note that the function $\xi$ satisfies $\xi (p,\beta) = - \phi \big(
1/p -1 , \beta\big)$.
\begin{theorem}\label{thmfive}
There exists $\delta' > 0$ and $p_0 \in (0,p_c)$ such that $p \in
(p_0,p_c)$ implies that
\begin{displaymath}
\lim_{n \to \infty}{\frac{- \log \mathbb{P}_p(\vert C(0) \vert = n)}{n}}
\end{displaymath}
exists and is given by 
\begin{displaymath}
\inf_{\beta \in (\lapha,\lapha + \delta')}{- \log g(\beta) + \xi(p,\beta)}.
\end{displaymath}
\end{theorem}
{\bf Proof}
We may write $\mathbb{P}_p (\vert C(0) \vert = n) = H_1 + H_2$, where
\begin{eqnarray}
 H_1 & = & \sum_{m=0}^{2(d-1)n -1}{\sigma_{n,m} p^n (1-p)^m}
 \nonumber \\
\textrm{and} \ H_2 & = & \sum_{m = 2(d-1)n}^{2(d-1)n+d}{\sigma_{n,m} p^n (1-p)^m}.
\end{eqnarray}
Note that, by Lemma \ref{lemmanex}, there exists $r \in (0,1)$ such that, for
all $p \in (0,1)$, $H_2 \leq (2d+1)r^n$.
To treat the quantity $H_1$, note that
\begin{displaymath}
H_1 = \sum_{m=0}^{2(d-1)n - 1}{a_n(m/n) \exp{-n(- \log g(m/n) + \xi(p,m/n))}}.
\end{displaymath}
We have that
\begin{displaymath}
H_1 \leq 2(d-1)L n^2 \exp{-n \gamma_p},
\end{displaymath}
where the quantity $\gamma_p$ is given by 
\begin{displaymath}
\gamma_p = \inf_{\beta \in [0,2(d-1))}{- \log g(\beta) + \xi(p,\beta)}.
\end{displaymath}
From Theorem \ref{thmtwo}, we see that the function
\mbox{$\beta \rightarrow - \log g(\beta) +
\xi(p,\beta)$} is continuous on $(0,2(d-1))$. Allied with the fact that for \mbox{$\beta \in (0,2(d-1))$},
\mbox{$\limsup{\frac{- \log a_n(\beta)}{n}} \leq 0$}, it follows that,
for $\epsilon > 0$, and for $n$ sufficiently large, $H_1 \geq
\exp{-n( \gamma_p + \epsilon)}$. Hence,
\begin{equation}\label{wtone}
\lim_{n \to \infty}{\frac{- \log H_1}{n}} = \gamma_p.
\end{equation} 
We now make the claim that there exists $p_0 \in (0,p_c)$ and $\delta'
> 0$ such that, for $p \in (p_0,p_c)$, $\gamma_p$ is given
by
\begin{equation}\label{pertin}
\gamma_p = \inf_{\beta \in [\lapha,\lapha+\delta']}{- \log
g(\beta) + \xi(p,\beta)}.
\end{equation}
Note that from
\begin{displaymath}
\frac{d}{d \beta} {\xi(p,\beta)} = \log (1 - 1/(1+\beta)) - \log(1-p),
\end{displaymath}
and $p_c = 1/(1 + \lapha)$, it follows that $\xi(p,\beta)$ is
decreasing in $\beta$ on $[0,\lapha]$, for $p \in (0,p_c)$. 
Let $\epsilon > 0$.
The upper semicontinuity of $g$ and the fact that $g(\beta) < 1$ for $\beta >
\lapha$ imply that we may find $\delta' > 0$ such that
\begin{displaymath}
\beta > \lapha + \delta' \ \textrm{implies that} \ g(\beta) < 1 -
\epsilon .
\end{displaymath}
Choose $\delta \in (0,\delta')$ such that
\begin{displaymath}
\beta \in (\lapha,\lapha + \delta) \ \textrm{implies that} \ g(\beta) > 1 -
\epsilon/2 .
\end{displaymath}
Let $p_0 = 1/(\lapha + \delta) - 1$. For $p \in (p_0,p_c)$,
\begin{equation}\label{wttwo}
- \log g(1/p - 1) + \xi(p,1/p -1) < - \log (1 - \epsilon/2),
\end{equation}
the second term on the left-hand-side being zero. 
For $\beta \in (\lapha + \delta' , 2(d-1)]$,
\begin{displaymath}
 - \log g(\beta) + \xi(p,\beta) \geq - \log (1 - \epsilon),
\end{displaymath}
since $\xi$ is non-negative. Hence, for $p \in (p_0,p_c)$,
(\ref{pertin}) holds, implying the claim.
Note that 
\begin{displaymath}
\liminf_{n \to \infty}{\frac{- \log H_2}{n}} \geq -
\log r,
\end{displaymath}
whereas, for $p \in (p_0,p_c)$, it follows from (\ref{wtone}) and (\ref{wttwo}) that 
\begin{displaymath}
\limsup_{n \to \infty}{\frac{- \log H_1}{n}} < - \log (1 -
\epsilon/2). 
\end{displaymath}
By choosing $\epsilon < 2(1-r)$, we obtain for such values of $p$,
\begin{displaymath}
\lim_{n \to \infty}{\frac{- \log \mathbb{P}_p(\vert C(0) \vert =
n)}{n}} = \inf_{\beta \in [\lapha,\lapha + \delta']}{- \log g(\beta) +
\xi(p,\beta)},
\end{displaymath} 
as required. $\Box$ \\
\subsection{Relating $\reta$ to the exponent for correlation size}
Theorem \ref{thmfive} allows us to deduce a scaling law that relates the combinatorially defined exponent $\reta$ to one which is defined directly from the percolation model.
\begin{definition}
Let $q:(0,p_c) \to [0,\infty)$ be given by
\begin{displaymath}
q(p) = \lim_{n \to \infty}{\frac{- \log \mathbb{P}_p(\vert C(0) \vert = n)}{n}}.
\end{displaymath}
\end{definition}
{\bf Remark} The existence of $q$ follows from a standard
subadditivity argument.  
\begin{definition} Define 
$\Omega_{+}^{\rrho} = \{ \gamma \geq 0 : \liminf_{p \uparrow p_c}  \frac{q(p)}{(p_c - p)^{\gamma}}  = \infty \}$ \ \textrm{and}\\
$\Omega_{-}^{\rrho} = \{ \gamma \geq 0 : \limsup_{p \uparrow p_c}{ \frac{q(p)}{(p_c - p)^{\gamma}}} = 0 \}$ 
\end{definition}
\begin{definition}
If $\sup{\Omega_{-}^{\rrho}} = \inf{\Omega_{+}^{\rrho}}$, then hypothesis $(\rrho)$ is said to hold, and $\rrho$ is defined to be equal to the common value.
\end{definition}
The function $q$ gives the exponential decay rate for the probability of observing a large cluster; this decay rate tends to zero as $p$ approaches $p_c$. Hypothesis $(\rrho)$ is introduced to describe how quickly that convergence occurs.

\begin{theorem}\label{thmsix}
Assume hypothesis $(\reta)$.
\begin{itemize}
\item Suppose that $\reta \in (1,2)$. Then hypothesis $(\rrho)$ holds and
$\rrho = 2$. 
\item Suppose that $\reta \in (2,\infty)$. Then hypothesis $(\rrho)$
holds and $\rrho = \reta$.
\end{itemize}
\end{theorem}
{\bf Proof}
Suppose that $\reta \in (1,2)$. Choose $\epsilon > 0$ so that $1 < \reta
- \epsilon < \reta + \epsilon < 2$. There exists constants $C_1,C_2>0$
such that, for $p \in (p_0,p_c)$ and $\beta \in (\lapha,\lapha +
\delta')$,
\begin{eqnarray}
 (\beta - \lapha)^{\reta + \epsilon} + C_1 \big( \beta - (1/p - 1)
 \big)^2 & \leq & - \log g(\beta) + \xi(p,\beta) \label{rtf} \\
 & \leq &  (\beta - \lapha)^{\reta - \epsilon} + C_2 \big( \beta - (1/p - 1)
 \big)^2. \nonumber
\end{eqnarray}
Applying Theorem \ref{thmfive}, we find that
\begin{equation}\label{rtfo}
 (\beta_p - \lapha)^{\reta + \epsilon} + C_1 \big( \beta_p - (1/p - 1)
 \big)^2  \leq q(p), 
\end{equation}
where $\beta_p \in [\lapha,\lapha +
\delta']$ denotes a value at which the infimum in the interval $[\lapha,\lapha 
\delta']$ of the first term in (\ref{rtf}) is  
attained. 
Let $y_p = 1/p - 1 - \lapha$, and let $\sigma_p$ satisfy $\beta_p = \lapha +
y_p^{\sigma_p}$. Then $\beta_p$ and $\sigma_p$ satisfy
\begin{eqnarray}
 (\reta + \epsilon) (\beta_p - \lapha)^{\reta + \epsilon - 1} & = & - 2 C_1
\big( \beta_p - (1/p - 1) \big) \nonumber \\
 (\reta + \epsilon) y_p^{\sigma_p ( \reta + \epsilon - 1)} & = &  2 C_1
\big( y_p - y_p^{\sigma_p} \big) \label{ghb}
\end{eqnarray}
Since $\beta_p \leq 1/p - 1$, $\sigma_p \geq 1$. From this and
(\ref{ghb}) follows $\liminf_{p \uparrow p_c}{\sigma_p} \geq
1/(\reta + \epsilon - 1)$. Applying (\ref{ghb}) again, we deduce that
 \mbox{$\lim_{p \uparrow p_c}{\sigma_p} =
1/(\reta + \epsilon - 1)$}.
Substituting $\sigma_p$ in (\ref{rtf}) yields
\begin{displaymath}
 y_p^{\sigma_p (\reta+ \epsilon)} + C_1 \big( y_p - y_p^{\sigma_p}
 \big)^2 \leq q(p).
\end{displaymath}
The facts that $\lim_{p \uparrow}{\sigma_p} > 1$ and $\lim_{p
\uparrow}{\sigma_p (\reta + \epsilon)} = (\reta + \epsilon)/(\reta +
\epsilon - 1) > 2$ imply that, for a small constant $c$, \mbox{$c (p_c
- p)^2 \leq q(p)$} for values of $p$ just less than $p_c$. 
 A similar
analysis in which $q(p)$ is bounded below by the infimum on the
interval $[\lapha,\lapha + \delta']$ of the third expression in
(\ref{rtf}) implies that for large $C$,  \mbox{$q(p) \leq C (p_c
- p)^2$}, in a similar range of values of $p$. Thus hypothesis $(\rrho)$ holds, and $\rrho = 2$. 

In the case where $\reta > 2$, let $\epsilon > 0$ be such that $\reta >
2 + \epsilon$. Defining $\sigma'_p$ by \mbox{$\beta_p = 1/p - 1 -
y_p^{\sigma'_p}$}, we find that 
\begin{equation}\label{qew}
 (\reta + \epsilon) 
\big( y_p - y_p^{\sigma'_p} \big)^{\reta + \epsilon - 1} =  2 C_1
 y_p^{\sigma'_p}.
\end{equation}
Note that $\beta_p \geq \lapha$ implies that $\sigma'_p \geq 1$. From
(\ref{qew}), it follows that $\liminf_{p \uparrow
p_c}{\sigma'_p} \geq \reta + \epsilon - 1$.  Since $\reta + \epsilon -
1 > 1$, applying (\ref{qew}) again shows that the limit $\lim_{p \uparrow
p_c}{\sigma'_p}$ exists and infact equals $\reta + \epsilon - 1$. Substituting $\sigma'_p$ in (\ref{rtf}) yields
\begin{displaymath}
 \big( y_p - y_p^{\sigma'_p}
 \big)^{\reta+ \epsilon} + C_1 y_p^{2 \sigma'_p} \leq q(p).
\end{displaymath}
The fact that $\liminf_{p \uparrow p_c}{\sigma'_p} > 1$ implies that
\mbox{$c (p_c - p)^{\reta + \epsilon} \leq q(p)$} for values of $p$
just less than $p_c$.  Making use of the
inequality $\reta > 2 + \epsilon$ in considering the infimum
of the third term appearing in (\ref{rtf}) yields in this case 
\mbox{$q(p) \leq C (p_c - p)^{\reta - \epsilon}$} for similar values of
$p$. Thus, since $\epsilon$
 may be chosen to be arbitrarily small, we find that, if $\reta > 2$, then hypothesis $(\rrho)$ holds, and
that $\rrho = \reta$. $\Box$ 

\section{Appendix: Proof of Theorem \ref{thmone}}
\subsection{Outline}
We begin with an outline of the argument for the existence of $f$. Note first of all that any lattice animal with $n$ edges has at most $2(d-1)n + 2d$ outlying edges. Let us restrict our attention to the case where $\beta \in (0,2(d-1))$. 

Take two animals with $n$ edges and surface-area-to-volume-ratio close
to $\beta$: $\gamma_1,\gamma_2 \in \Gamma'_{n,\lfloor \beta n \rfloor
}.$  
Translate $\gamma_2$ so that its lexicographically smallest vertex coincides with the lexicographically greatest vertex of $\gamma_1$, and form \mbox{$\gamma$ ($=\gamma_1 * \gamma_2$, say)}, by taking the union of $\gamma_1$ with this translate of $\gamma_2$. The animal $\gamma$ has $2n$ edges, and surface-area-to-volume ratio equal to $\beta$, with a correction, uniformly bounded in $n$, arising from an effect around the point of concatenation.
Since $\gamma_1$ and $\gamma_2$ can be identified from $\gamma$, we have demonstrated that
\begin{equation}\label{eqnafg}
\sigma'_{2n,\lfloor 2n\beta \rfloor} \geq (\sigma'_{n,\lfloor n\beta \rfloor})^2,
\end{equation}
provided that we overlook the small correction just mentioned.
Equation \ref{eqnafg} implies that
\begin{displaymath}
f'_{2n}(\beta) \geq f'_n(\beta),
\end{displaymath}
from which the existence of $f$ follows, given the fact that there is an exponential growth rate for the total number of lattice animals up to translation, as a function of size.

How to transform this argument into a proof?

\begin{lemma}\label{lemmadt}
Let $\{ a_n: n \in \mathbb{N} \}$ be a sequence of positive numbers for which there exists $k \in \mathbb{N}$ such that for all $n,m \geq 1$, we have $ a_{n+m+k} \geq a_{n}a_{m} $ and $ \sup{a_n^{1/n}} < \infty $.\\
Then $\lim_{n \to \infty}{a_n^{1/n}}$ exists, and
\begin{displaymath}
a_{m-k}^{1/m} \leq \lim_{n \to \infty}{{a_n}^{1/n}} \quad \textrm{for all m.}
\end{displaymath}
\end{lemma}
{\bf Remark} This lemma and its proof appear in \cite{MR95m:82076}.

The context in which Lemma \ref{lemmadt} will be applied is that of a fixed $\beta \in (0,2(d-1))$. Setting $a_n = \sigma'_{n,\lfloor \beta n \rfloor}$, the hypothesis of the lemma may be rewritten:\\
there exists $K \in \mathbb{N}$, such that, for all $n,m \in \mathbb{N}$,
\begin{displaymath}
\sigma'_{n+m+K,\lfloor \beta (n+m+K) \rfloor }\geq \sigma'_{n,\lfloor \beta n \rfloor } \sigma'_{m,\lfloor \beta m \rfloor}.
\end{displaymath}

From this form, it is possible to explain the technique by which the concatenation argument may be made rigorous. The animal $\gamma$ formed by concatenating $\gamma_1$ and $\gamma_2$ may have a surface-area-to-volume ratio which needs to be adjusted. The adjustment is effected by the further concatenation of a specific construction, carefully designed to manipulate the surface area by the small amount required. It is the additional edges in this construction that contribute the extra $K$ that appears in the statement of Lemma \ref{lemmadt}.

Having explained the ideas, we introduce some of the objects required for the proof of Theorem \ref{thmone}.

\subsection{Constructions}

Throughout this section, we let $d \in \mathbb{N}$, with  $d \geq 2$. Various constructions will be required.

\begin{table}
\begin{tabular}{|r|l|}
\hline
$\rho^{a_1, \ldots ,a_k}$ & $(a_1, \ldots , a_k) \in \{0,1\}^{\{1,\ldots,k\}}$ \quad $ k \in \mathbb{N} $ \\
$\rho^{a_1, \ldots ,a_k}_d$ & $(a_1, \ldots , a_k) \in \{0,1\}^{\{1,\ldots,k\}}$ \quad $ k \in \mathbb{N} $ \\
$ S_k^i $ & $ k,i \in \mathbb{N}$\quad$  i \in \{0, \ldots, k-2 \}$ \\
$\phi_{q_1,\ldots,q_d} $ & $ (q_1, \ldots , q_d) \in \mathbb{N}^d $ \\
\hline
\end{tabular}
\caption{The constructions and the parameters that index them}
\end{table}

The first and third of these constructions are planar, and are embedded in the $d$-dimensional lattice by the injection $\mathbb{Z}^2 \subseteq \mathbb{Z}^d$ : $ (x_1,x_2)  \to (x_1,x_2,0, \ldots,0)$. The second construction is formed from the first by adding a collection of edges that do not lie in the embedded plane.

Each of the four constructions is now defined.

\subsubsection{$\rho^{a_1, \ldots ,a_k}$}

\begin{figure}
\begin{center}
\includegraphics{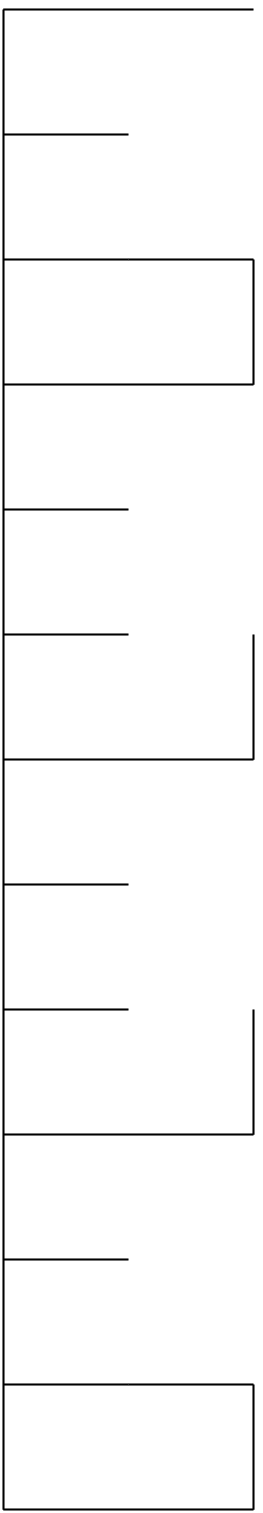}
Construction $\rho^{0,1,1,0}$.
\end{center}
\end{figure}

Let $k \in \mathbb{N}$ and $\{a_1, \ldots, a_k \} \in \{0,1\}^{\{1,\ldots,k\}}$.
This construction is planar. It contains the edges
\begin{displaymath}
\begin{array}{lll}
(0,0) - (1,0), & (1,0)-(2,0), & (0,0)-(0,1),   \\
(2,0) - (2,1), & (0,1)-(1,1), & (0,1)-(0,2), \\ 
(0,2) - (0,3). & &
\end{array} 
\end{displaymath}
It also contains the translates of these edges by $(0,3),(0,6),\ldots,(0,3(k-1))$. Finally, the edges $(0,3k)-(1,3k)$ and $(1,3k)-(2,3k)$ are present.

That is, the construction consists a collection of $k$ chambers in the plane. The seven edges listed above lie in the first chamber; their translates by $(0,3(j-1))$ lie in the $j^{\rm th}$ chamber, for $j \in \{ 1, \ldots, k \}$. The final two edges lie in the $k^{\rm th}$ chamber. In this way, the chambers are ordered $\{ 1,\ldots,k \}$ in increasing $y$-coordinate. Into each chamber, an edge is added. The position of that edge is determined by the status of $a_i$.

If $a_1=0$, then edge $(1,1)-(2,1)$ is added.
If $a_1=1$, then edge $(0,2)-(1,2)$ is added.

Similarly, in the other cases: if $a_j=0$, then edge $(1,1 + 3(j-1)),(2,1 + 3(j-1))$ is added.
If $a_j=1$, then edge $(0,2 + 3(j-1))-(1,2+ 3(j-1))$ is added.

In the figure, $\rho^{0,1,1,0}$ is depicted.

\subsubsection{$\rho^{a_1, \ldots ,a_k}_d$}

Once again, let  $k \in \mathbb{N}$ and $\{a_1, \ldots, a_k \} \in \{0,1\}^{\{1,\ldots,k\}}$.

Having embedded the first construction $\rho^{a_1, \ldots ,a_k}$ in $\mathbb{Z}^d$, some extra edges are added. These will amount to $4(d-2)k$ in total; they will be divided into $k$ classes, corresponding to each chamber of $\rho$,  with $4(d-2)$ in each class.

Class $1$ consists of the following edges:

\begin{displaymath}
\begin{array}{ccc}
(0,2,0,\ldots,0) & - & (0,2,1,\ldots,0) \\
(0,2,0,\ldots,0) & - & (0,2,0,1,\ldots,0)\\
 \vdots          &   &                 \\ 
(0,2,0,\ldots,0) & - & (0,2,0,\ldots,1)
\end{array} 
\end{displaymath}

and:
\begin{displaymath}
\begin{array}{ccc}
(0,2,1,\ldots,0)   & - & (1,2,1,0,\ldots,0)\\
(0,2,0,1,\ldots,0) & - & (1,2,0,1,\ldots,0)\\
\vdots             &   &                  \\
(0,2,0,\ldots,1)   & - & (1,2,0,\ldots,0,1)\\
\end{array} 
\end{displaymath}

and:
\begin{displaymath}
\begin{array}{ccc}
(0,2,0,\ldots,0) & - & (0,2,-1,\ldots,0)\\
(0,2,0,\ldots,0) & - & (0,2,0,-1,\ldots,0)\\
\vdots           &   &                   \\
(0,2,0,\ldots,0) & - & (0,2,0,\ldots,-1)\\
\end{array} 
\end{displaymath}

and:
\begin{displaymath}
\begin{array}{ccc}
(0,2,-1,0,\ldots,0)  & - &(1,2,-1,0,\ldots,0)\\
(0,2,0,-1,0,\ldots,0)& - & (1,2,0,-1,0,\ldots,0)\\
\vdots            &   &               \\
(0,2,0,\ldots,0,-1)  & - & (1,2,0,\ldots,0,-1)\\
\end{array} 
\end{displaymath}

Class $j$ consists of the translate Class $1$ $ + \, \, (0,3(j-1),0,\ldots,0)$.

\subsubsection{$S_k^i$}

Let $k \in \mathbb{N}$ and $i \in \{0,\ldots, k-2 \}$.

$S_k^i$ contains a frame of edges, whose horizontal elements are given by:

\begin{displaymath}
\begin{array}{l}
(0,0)-(1,0),(1,0)-(2,0),\ldots,(k-1,0)-(k,0),\\
(0,2)-(1,2),(1,2)-(2,2),\ldots,(k-1,2)-(k,2),\\
\end{array}
\end{displaymath}
along with the four vertical edges $(0,0)-(0,1)$, $(0,1)-(0,2)$, $(k,0)-(k,1)$ and $(k,1)-(k,2)$.

It also contains the edges:
$ \Big\{ (1,0)-(1,1),\ldots,(i,0)-(i,1) \Big\} $.

\subsubsection{$\phi_{q_1,\ldots,q_d}$}

Let $(q_1, \ldots , q_d) \in \mathbb{N}^d$.

The vertex set $V$ of $\phi$ comprises $\{0,1,\ldots,q_1\} \times \{0,\ldots,q_2\} \times \ldots \times \{0,\ldots,q_d\}$. The edge set consists of all those edges both of whose endpoints lie in $V$.

\subsection{Calculating statistics of constructions}

The constructions will be concatenated for appropriate choices of the parameters that index them to form a structure which will be added to two joined lattice animals in a way that adjusts surface area correctly. As such, we must calculate the number of edges and the number of outlying edges of the constructions. The notation $e$ and $o$ will be used to denote these quantities respectively.
\begin{itemize}
\item{$\rho^{a_1, \ldots ,a_k}$}

We have that the number of edges $e(\rho^{\vec a})$ of the construction $\rho^{\vec a}$ is given by
\begin{displaymath}
e(\rho^{\vec a}) = 8k + 2,
\end{displaymath}
and that its number of outlying edges is equal to
\begin{displaymath}
o^2(\rho^{\vec a})= 11k + \sum_{i=1}^{k}{a_i} + 8 .
\end{displaymath}
The notation $o^2$ is used to emphasise that we are considering the
construction $\rho^{\vec a}$ as a subset of $\mathbb{Z}^2$, rather
than its embedding in $\mathbb{Z}^d$. 
\item{$\rho^{a_1, \ldots ,a_k}_d$}

In this case, the number of edges is given by
\begin{eqnarray}
 e(\rho^{\vec a}_d) & = & e(\rho^{\vec a}) + 4(d-2)k {} \nonumber \\
      & = & 8k + 2 + 4(d-2)k . \nonumber
\end{eqnarray}

The outlying edges $O(\rho^{\vec a}_d)$ of $\rho^{\vec a}_d$ can be divided as follows:
\begin{displaymath}
O(\rho^{\vec a}_d) = O_1 \cup O_2
\end{displaymath}
where
$O_2 = \{ e \in O(\rho^{\vec a}_d):$ $e$ does not lie in the $(x_1,x_2)$-plane and is incident to an edge of $\rho^{\vec a}_d$ not lying in the $(x_1,x_2)$-plane $\}$
and $O_1 = O(\rho^{\vec a}_d) - O_2$.

As such,
\begin{eqnarray}
\vert O_1 \vert & = & o^2(\rho^{\vec a}) + 2(d-2)[V(\rho^{\vec a}) - \sum_{i=1}^{k}{a_i} - k] , \nonumber \\
\textrm{and} \ \vert O_2 \vert & = & 2(d-2)(4d-3)k . \nonumber
\end{eqnarray}
So,
\begin{displaymath}
o(\rho^{\vec a}_{d}) = 11k + \sum_{i=1}^{k}{a_i} + 8 + 2(d-2)(6k+3) + 2(d-2)(4d-3)k .
\end{displaymath}
\item{$S_k^i$}

We have that
\begin{eqnarray}
e(S_k^i) & = & 2k + 4 + i \nonumber \\ 
o^2 (S_k^i) & = & 4k + 8  \nonumber \\
o(S_k^i) & = & 4k + 8 + 2(d-2)(2(k+2)+i) \nonumber
\end{eqnarray}
\item{$\phi_{q_1,\ldots,q_d}$}

To calculate the number of edges and outlying edges in this construction, we divide its vertices into classes, as follows.

\begin{displaymath}
\{ 0,1,\ldots,q_1 \} \times \ldots \times \{0,1,\ldots,q_d \} = \cup_{0 \leq k \leq d} \cup_{j \in \{1,\ldots,{{d \choose k}}\}}C(k,j)\\
\end{displaymath}
where $C(k,j)$ consists of those elements $ \{c_1, \ldots, c_d \} $ of the vertex set such that $T_j := \{ i \in \{ 1, \ldots, d \}: c_i = q_i \} $ is equal to the $j^{\rm th}$ lowest $k$-subset of $\{ 1,\ldots,d \} $ in the lexicographical ordering.

To each vertex $x$ in $C(k,j)$ (for $0 \leq k \leq d-1 , j \in \{ 1, \ldots , {d \choose k} \} $), write $\sigma(x)$ for the set of $d-k$ edges emanating from $x$ towards a lexicographically greater endpoint which lies in $V(\phi_{q_1,\ldots,q_d})$.

Then $e(\phi)$, the edge set of the construction, can be represented as a disjoint union of the sets $\{ \sigma(x): x \in V(\phi) \}$.

As such,
\begin{displaymath}
e(\phi_{q_1,\ldots,q_d})= \sum_{k=0}^{d-1}{\sum_{j=1}^{d \choose k}{\vert C(k,j) \vert (d-k)}}.
\end{displaymath}

For $j \in \{ 1, \ldots, {d \choose k} \}$, $\vert C(k,j) \vert = \prod_{i \in \{ 1, \ldots, d \} - T_j } {q_i}$.

Let $S_j$ denote the $j^{\rm th}$ symmetric polynomial in $\{ q_1, \ldots, q_k \} $. The notation is displayed in Table $2.2$. 
\\
\begin{table}\label{sympoly}
\begin{center}
\begin{tabular}{l|l}
\hline
$ S_0 $ & $ q_1 \times \ldots \times q_d $ \\
$ S_1 $ & $ \sum_{i=1}^{d}{\prod_{j \neq i}{q_j}} $ \\
$ \ldots $ \\
$ S_{d-1} $ & $ q_1 + \ldots + q_d $ \\
\hline
\end{tabular}
\caption{Formulae for $S_j(q_1,\ldots,q_d)$}
\end{center}
\end{table}

We have that
\begin{displaymath}
e( \phi_{q_1, \ldots, q_d}) = \sum_{k=0}^{d-1}{(d-k)S_k(q_1,\ldots,q_d)}.
\end{displaymath}

The outlying edge set of $\phi_{q_1,\ldots,q_d}$, $O(\phi)$, can be identified with the collection $G$ of endpoints of its members that do not lie in $V(\phi)$.
 The set $G$ may be decomposed as 
\begin{displaymath}
G = \cup_{k=1}^{d}{G_k \cup G^k} ,\\
\end{displaymath}
where
\begin{eqnarray}
G_k & = & \{ (a_1,\ldots,a_d): \textrm{for} \, j \neq k, a_j \in \{ 0,\ldots, q_j \} , a_k = -1 \} \nonumber \\
\textrm{and} \, \, G^k & = & \{ (a_1,\ldots,a_d): \textrm{for} \, j \neq k, a_j \in \{ 0,\ldots, q_j \} , a_k = q_k + 1 \} \nonumber .
\end{eqnarray}
Now, 
\begin{displaymath}
\vert G_k \vert = \vert G^k \vert = \prod_{i=1,i \neq k}^{d}{(q_i + 1)},
\end{displaymath}
which implies that 
\begin{displaymath}
o(\phi_{q_1,\ldots,q_d}) = \vert G \vert = 2 \sum_{k=1}^{d}{\prod_{i=1.i \neq k}^{d}{(q_i + 1)}}.
\end{displaymath}
\end{itemize}

\begin{table}\label{tablesumm}
\begin{tabular}{l|l|l}
\hline
E & e(E) & o(E)  \\
\hline
$\rho^{a_1,\ldots,a_k}_{d}$ & $8k + 2 + 4(d-2)k $ & $11k + \sum_{i=1}^{k}{a_i} + 8 + 2(d-2)[(4d+3)k+3] $\\
$S_k^i$ & $2k+4+i$ & $2(d-2)(2k+4+i) + 8 + 4k$ \\
$\phi_{q_1, \ldots,q_d}$ & $\sum_{k=0}^{d-1}{(d-k)S_k(q_1,\ldots,q_d)}$ & $2 \sum_{k=1}^{d} {\prod_{i=1,i \neq k}^{d}{(q_i + 1)}} $\\
\hline
\end{tabular}
\caption{Summary of statistics of constructions}
\end{table}

\subsection{Operation *}
At the heart of the argument that underlies the proof of Theorem \ref{thmone} is the idea of joining together two large lattice animals of equal surface-area-to-volume ratio to obtain a new animal, of twice the size, and of very similar surface-area-to-volume ratio as its two constituents. In this section, we define the concatenation procedure that will be employed, and calculate the precise surface-area-to-volume ratio of the animal that results from using it.
\begin{definition}
Let $\Gamma'$ denote the collection of lattice animals whose
lexicographically minimal vertex is the origin.
Let $*:\Gamma' \times \Gamma' \to \Gamma'$ be defined as follows. Let $\gamma_1,\gamma_2 \in \Gamma'$. Let $\gamma = \gamma_1 \cup (\gamma_2 + v)$, where $(\gamma_2 + v)$ is the translate of $\gamma_2$ whose lexicographicxally lowest vertex is equal to the lexicographically greatest vertex of $\gamma_1$. Then $\gamma = \gamma_1 * \gamma_2$.
\end{definition}

\begin{definition}
Let $R_2$ be the two-edged lattice animal whose edges are \\
$(0,\ldots,0)-(1,0,\ldots,0),(1,0,\ldots,0)-(2,0,\ldots,0)$.
\end{definition} 

\begin{lemma}\label{littlelem}
Let $\gamma_1 \in \Gamma'_{n_1,m_1}, \gamma_2 \in \Gamma'_{n_2,m_2}$.  If $\gamma_1 * R_2 * \gamma_2 \in \Gamma'_{n_3,m_3}$, then $ n_3 $ and $ m_3 $ satisfy:\\
$n_3 = n_1 + n_2 + 2$\\
$m_3 = m_1 + m_2 + 2(d-2)$.
\end{lemma}

{\bf Proof}
The formula for $n_3$ is trivial. The additional term of $2(d-2)$ in the formula for $m_3$ is explained by the fact that $2(d-1)$ extra outlying edges arise emanating from the midpoint of $R_2$; the edge of $R_2$ that touches the origin blocks an outlying edge of $\gamma_1$, and the other edge of $R_2$ blocks an outlying edge of $\gamma_2$. 

\subsection{Assembling the components}
As already mentioned, the idea
for proving the existence of $f$ is to add to the join of two large lattice animals an additional structure designed to correct the slight error in surface-area-to-volume ratio that arises from the joining. Having described the constructions which form the building blocks of this structure, calculated their relevant properties, and defined the concatenation procedure, we are now in a position to describe the form of the structure we will use more explicitly.

Firstly, we fix surface-area-to-volume ratio at some value $\beta \in (0,2(d-1))$. Choosing $n_1, n_2 \in \mathbb{N}$ and writing $m_1 = \lfloor \beta n_1 \rfloor$ and $m_2 = \lfloor \beta n_2 \rfloor$, we let $\gamma_1 \in \Gamma'_{n_1,m_1}, \gamma_2 \in \Gamma'_{n_2,m_2}$. In introducing the corrected form of the joined animals, we make our first choice for the parameters, by insisting that $2d-1$ quantities $a_i$ will be used in the construction $\rho_d^{\vec a}$. That is, we write
\begin{displaymath}
\tau = \gamma_1 * R_2 * \gamma_2 * R_2 * \rho_d^{\{a_1,\ldots,a_{2d-1}\}} * R_2 * S_k^i * R_2 * \phi_{\{q_1,\ldots,q_d \}}.
\end{displaymath}

The correction structure appearing in $\tau$ comprises three parts, called $\rho$, S, and $\phi$. This new structure will contribute new edges, and will affect the number of outlying edges. Precisely how depends on the various parameters which define the three building blocks, and we perform the necessary calculations shortly. 
Let's explain the role of each of these building blocks in determining the required correction construction. 
We have a fixed $\beta$; in adding a structure with $K$ edges, we must ensure that the structure also contributes roughly $\beta K$ outlying edges. In addition to this, the structure must manipulate the number of outlying edges to precisely the required value.\\
The role of $\phi$ and $S$ is to ensure that the additional structure has roughly the right surface-area-to-volume ratio. A choice of the parameters that index $\phi$ will shortly be made so that this construction takes the form of a cube; the surface-area-to-volume ratio of this cube may be made arbitrarily small by taking the side-length high enough. The construction $S$, on the other hand, is an animal with a high surface-area-to-volume ratio. By allowing $\phi$ and $S$ to have the correct relative sizes, we ensure that the additional structure has roughly the right surface-area-to-volume ratio. Having done this, the parameter choice $i$ for the construction $S$ is made to manipulate the number of outlying edges of the whole body to within a finite discrepancy of the required value. 
The most delicate adjustment is made by using $\rho$: by choosing the values of $a_i$ correctly, we will adjust the surface-area-to-volume ratio to the precise value required. The construction $\rho$ was designed with this purpose in mind: changing an $a_i$ from 0 to 1 effects an increase in the number of outlying edges of precisely one. The overhanging structures that lie outside of the plane $\mathbb{Z}^2 \subseteq \mathbb{Z}^d$ were put in place to ensure that this property holds. 

Let us compute $e(\tau)$ and $o(\tau)$ for the case of general parameter values, before beginning to specify the values that they should take for given $\beta$.

\begin{eqnarray}
e(\tau) & = & e(\gamma_1) + e(\gamma_2) {}\label{ghjj} \\
 & & {} + e(\rho_d^{\{a_1,\ldots,a_{2d-1}\}}) + e(S_k^i) + e(\phi_{\{q_1,\ldots,q_d\}}) + 8 {} \nonumber \\
o(\tau) & = & o(\gamma_1) + o(\gamma_2) {}\label{ghjk}\\
 & & {} + o(\rho_d^{\{a_1,\ldots,a_{2d-1}\}}) + o(S_k^i) +
 o(\phi_{\{q_1,\ldots,q_d\}}) + 8(d-2) , \nonumber
\end{eqnarray}
where the final terms of $8$ and $8(d-2)$ in these expressions occur
as a result of the presence of the four joining structures $R_2$ in the
construction $\tau$; we are using Lemma \ref{littlelem} to compute
these terms.
The number of edges and outlying edges in each of the constructions
$\rho_d^{\{a_1,\ldots,a_{2d-1}\}}, S_k^i$ and
$\phi_{\{q_1,\ldots,q_d\}}$ have been computed; the resulting
expressions are shown in Table $2.3$.
We now substitute these expressions in the two equations, \ref{ghjj}
and \ref{ghjk}. Recalling that $e(\gamma_1)=n_1$, $e(\gamma_2)=n_2$,
$o(\gamma_1)=m_1$, $o(\gamma_2)=n_2$, we obtain,   
\begin{displaymath}
e(\tau)  =  n_1 + n_2 + 4d(2d-1) + 14 + 2k  + i + \sum_{j=0}^{d-1}{(d-j)S_j(q_1,\ldots,q_d)},
\end{displaymath}
and
\begin{eqnarray}
o(\tau) & = & m_1 + m_2 + 11(2d-1) + \sum_{j=1}^{2d - 1}{a_j} + 8 +
2(d-2)[(4d+3)(2d-1) + 3] \nonumber \\
 & & {} + 2(d-2)(2k + 4 + i) + 8 + 4k + 2\sum_{j=1}^{d}{ \prod_{r \in
 \{ 1, \ldots ,d \}, r \neq j}{(q_r + 1) }} + 8(d-2) \nonumber \\ 
   & = & m_1 + m_2 + \sum_{j=1}^{2d-1}{a_j} + 2(d-2)[(4d+3)(2d-1) + 11] + 16{}\nonumber \\
 & + & 4(d-1)k + 2(d-2)i + 2\sum_{j=1}^{d}{ \prod_{r \in \{ 1, \ldots ,d \}, r \neq j}{(q_r + 1) }} + 11(2d-1) . \nonumber 
\end{eqnarray}   
Our aim is to choose the parameters in such a way that there exists $K$ satisfying
\begin{displaymath}
(e(\tau),o(\tau))= (n_1 + n_2 + K, \lfloor \beta (n_1 + n_2 + K) \rfloor ).
\end{displaymath}
We now set $q_j = \lfloor (ak)^{1/d} \rfloor$ for all $j \in \{ 1, \ldots, d \}$, with $k \in \mathbb{N}$, \mbox{$i \in \{ 0, \ldots, k-1 \}$} and $a \in \mathbb{R}$. It is this choice that determines that the construction $\phi$ takes the form of a $d$-cube. The value of $a \in \mathbb{R}^{+}$ and of $i \in \{ 1,\ldots, k-2 \}$ will be determined shortly.   
We require that $K$ satisfies
\begin{equation}\label{eqnk}
K = 4d(2d-1) + 14 + 2k + i + \sum_{j=0}^{d-1}{(d-j)S_j(\lfloor (ak)^{1/d} \rfloor,\ldots,\lfloor (ak)^{1/d} \rfloor)}
\end{equation}
and
\begin{eqnarray}\label{eqns}
 & & \lfloor \beta (n_1 + n_2 + K) \rfloor - \lfloor \beta n_1 \rfloor - \lfloor \beta n_2 \rfloor \\ 
 & = & {} \sum_{j=1}^{2d-1}{a_j} + 2(d-2)\{(4d+3)(2d-1) +11\} + 16 + 4(d-1)k {}\nonumber \\ 
 & + & 2(d-2)i + 2 \sum_{j=1}^{d} { \prod_{r \in \{ 1, \ldots, d \}, r \neq j } {\{\lfloor (ak)^{1/d} \rfloor + 1 \}}} + 11(2d-1) \nonumber,
\end{eqnarray}
where we have used $m_1 = \lfloor \beta n_1 \rfloor$ and $m_2 = \lfloor \beta n_2 \rfloor$.
Note that 
\begin{displaymath}
\lfloor \beta (n_1 + n_2 + K) \rfloor - \lfloor \beta n_1 \rfloor  - \lfloor \beta n_2 \rfloor  - \lfloor \beta K \rfloor \in \{ 0,1,2 \}.
\end{displaymath}
Set $a_1 + a_2$ equal to this expression.
Doing so reduces the condition stated in \ref{eqns} to
\begin{eqnarray}
\lfloor \beta K \rfloor & = & \sum_{i=3}^{2d-1}{a_i} + 2(d-2)\{(4d+3)(2d-1) +11\}+ 16 + 4(d-1)k {}\nonumber\\
  & & {} + 2(d-2)i + 2d(\lfloor (ak)^{1/d} \rfloor +1)^{d-1} + 11(2d-1) . \nonumber
\end{eqnarray}
Note that Equation \ref{eqnk} is given by
\begin{displaymath}
K = 4d(2d-1) + 14 + 2k + i + d(\lfloor (ak)^{1/d} \rfloor)^{d} + \overline{S}(ak) , 
\end{displaymath}
where
\begin{displaymath}
\overline{S}(ak)= \sum_{j=1}^{d-1}{(d-j)S_j(\lfloor (ak)^{1/d} \rfloor,\ldots,\lfloor (ak)^{1/d} \rfloor)} .
\end{displaymath}
We proceed to fix the parameter $a \in \mathbb{R}^{+}$ that determines the side-length of the $d$-cube $\phi_{\lfloor (ak)^{1/d} \rfloor, \ldots, \lfloor (ak)^{1/d} \rfloor}$.
\begin{itemize}
\item if $\beta \in (0,2(d-2))$, then let $a$ be a positive value lying in the interval
\begin{displaymath}
\Big( \frac{4(d-1) - 2 \beta}{\beta d}, \frac{4(d-1) + 2(d-2) - 3 \beta}{\beta d} \Big); 
\end{displaymath}
\item if $\beta \in (2(d-2),2(d-1))$, then let $a$ be a positive value  lying in the interval
\begin{displaymath}
\Big( \frac{4(d-1) + 2(d-2) -3 \beta}{\beta d}, \frac{4(d-1) - 2 \beta}{\beta d} \Big). 
\end{displaymath}
\end{itemize}
The case where $\beta$ takes the value $2(d-2)$ will be discussed later.
\begin{lemma}\label{lemgyu}
Let $\beta \in (0,2(d-1))-\{2(d-2)\}$. For $k \in \mathbb{N}$, let $I$ denote the interval given by
\begin{displaymath}
 I  = \left\{ \begin{array}{ll} \{1,\ldots,\lfloor k(\beta - 2(d-2)) \rfloor \} & \textrm{if  $\beta \in (2(d-2),2(d-1))$}\\
  \{ \lfloor k(\beta - 2(d-2)) \rfloor + 1, \ldots, -1 \} & \textrm{if  $\beta \in (0,2(d-2))$}
\end{array} \right.
\end{displaymath}
There exists $k_0$ such that, for all $k > k_0$, the expression
\begin{eqnarray}
& & 4(d-1)k  +  2d(\lfloor (ak)^{1/d} \rfloor +1)^{d-1} {} \nonumber \\
& + & 2(d-2)[(4d+3)(2d-1)+11] + 16 + 11(2d-1) {} \nonumber\\
& - & \Big\lfloor \beta \Big( 4d(2d-1) + 2k + 14 + \sum_{j=0}^{d-1}{(d-j)S_j \big( \lfloor (ak)^{1/d} \rfloor,\ldots,\lfloor (ak)^{1/d} \rfloor \big)} \Big) \Big\rfloor \nonumber
\end{eqnarray}
lies in the interval $I$.
\end{lemma}
{\bf Proof}
Let $T$ denote the expression appearing in the statement of the lemma.
Then
\begin{displaymath}
T = k(4(d-1) - (2+da)\beta ) + R,
\end{displaymath}
where $R$ is $O(k^{1-(1/d)})$ as $k \to \infty$.
The condition on $a$ implies that 
$[4(d-1) - (2+da) \beta]k$ lies in the interval $I$.
Hence, $T \in I$ for $k$ sufficiently large, as required. $\Box$ \\
We choose $k$ to be any integer greater than $k_0$, the constant of Lemma \ref{lemgyu}. In preparing to set the value of the parameter $i$, note that, by substituting the expression for $K$ that appears in Equation \ref{eqnk},
\begin{eqnarray}
& & \sum_{j=3}^{2d-1}{a_j} + 2(d-2)\{(4d+3)(2d-1) +11\} + 16 + {4(d-1)}k + 2(d-2)i  {} \nonumber \\
& + & 2d(\lfloor (ak)^{1/d} \rfloor +1)^{d} + 11(2d-1) - \lfloor \beta K \rfloor  {} \nonumber \\
& = & {} \sum_{i=3}^{2d-1}{a_i} + 2(d-2)\{(4d+3)(2d-1) +11\} + 16 {} \nonumber \\
& + & {4(d-1)}k + 2(d-2)i + 2d(\lfloor (ak)^{1/d} \rfloor +1)^{d} + 11(2d-1) - \lfloor i \beta \rfloor - \delta  {} \nonumber \\
& - & \Big\lfloor \beta \Big( 4d(2d-1) + 14 + 2k + \sum_{j=0}^{d-1}{ (d-j)S_j(\lfloor (ak)^{1/d} \rfloor,\ldots,\lfloor (ak)^{1/d} \rfloor ) } \Big) \Big\rfloor {} \nonumber ,
\end{eqnarray}
where $ \delta \in \{0,1\}$.
Set $a_3 = \delta$.

\begin{lemma}
Let $\beta \in (0,2(d-1)) - \{2(d-2) \}$. Let $a \in \mathbb{R}^{+}$, $a_1, a_2, a_3 \in \{ 0,1 \}$, and $k > k_0$ be as specified. Then we may choose $i \in \{ 1, \ldots, k-2 \}$, and $\{a_4,\ldots,a_{2d-1} \} \in {( 0,1 )}^{2d-4}$ such that
\begin{eqnarray} 
 & & \lfloor i \beta \rfloor - 2(d-2)i \, =  \sum_{j=4}^{2d-1}{a_j} + 2(d-2)[(4d+3)(2d-1)+11] {}\nonumber\\
 & + & 16 + 4(d-1)k + 2d({\lfloor (ak)^{1/d} \rfloor +1})^{d-1} + 11(2d-1) {}\nonumber \\
 & - & \, \Big\lfloor \beta \Big( 4d(2d-1) + 2k + 14 + \sum_{j=0}^{d-1}{(d-j)S_j(\lfloor (ak)^{1/d} \rfloor,\ldots,\lfloor (ak)^{1/d} \rfloor)} \Big) \Big\rfloor \nonumber 
\end{eqnarray}
\end{lemma}
{\bf Proof}
Consider the function $\{1,\ldots,k-2\} \to \mathbb{Z}$ given by
\begin{displaymath}
i \to (\lfloor i \beta \rfloor - 2(d-2)i ).
\end{displaymath}
Note that every element of the interval $I$, as defined in Lemma \ref{lemgyu}, is at most $d-2$ from some element of the image of this map. \\
We know by Lemma \ref{lemgyu} that $Q \in I$, where the expression $Q$ is given by 
\begin{eqnarray}
 & & 2(d-2)[(4d+3)(2d-1)+11]+16 {}\nonumber \\
 & + & 4(d-1)k + 2d({(\lfloor(ak)^{1/d}\rfloor+1)}^{d-1}) + 11(2d-1) {}\nonumber\\
 & - & \Big\lfloor \beta \Big( 4d(2d-1) + 14 + 2k + \sum_{j=0}^{d-1}{(d-j)S_j(\lfloor(ak)^{1/d}\rfloor,\ldots,\lfloor(ak)^{1/d}\rfloor)} \Big) \Big\rfloor . \nonumber 
\end{eqnarray}
So, we may choose $i \in \{1,\ldots,k-2\}$ such that

\begin{displaymath}
Q  \leq \lfloor i \beta \rfloor - 2(d-2)i \leq Q + 2(d-2) . 
\end{displaymath}
Set $ \sum_{j=4}^{2d-1}{a_i} = \lfloor i \beta \rfloor - 2(d-2)i - Q$, for this choice of $i$, to prove the result. $\Box$

To summarise, we have chosen the parameters that index each of the building blocks that define the correction construction that is added to \mbox{$\gamma_1 * R_2 * \gamma_2$} in forming $\tau$.

We let $k = k_0 +1$, and set  

\begin{displaymath}
K = 4d(2d-1) + 14 + 2k + i + \sum_{j=0}^{d-1}{(d-j)S_j(\lfloor (ak)^{1/d} \rfloor,\ldots,\lfloor (ak)^{1/d} \rfloor)},
\end{displaymath}
where $a \in \mathbb{R}$ and $i \in \{0, \ldots, k-1\}$. 
With these choices, the map 
\begin{displaymath}
\Gamma'_{n_1,\lfloor \beta n_1 \rfloor}*\Gamma'_{n_2,\lfloor \beta n_2 \rfloor} \to \Gamma'_{n_1 + n_2 + K,\lfloor \beta (n_1 + n_2 + K) \rfloor}: (\gamma_1,\gamma_2) \to \tau
\end{displaymath}
is a well-defined injection.

It is time to apply Lemma \ref{lemmadt}. Setting $a_n =
\sigma'_{n,\lfloor \beta n \rfloor}$, the hypothesis of the lemma is
satisfied. Defining temporarily $f(\beta)$ as $\lim_{n \to
\infty}{f'_n(\beta)}$, the conclusion of the lemma may be written,
\begin{displaymath}
f(\beta) \ \textrm{exists for} \ \beta \in (0,2(d-1)) - \{2(d-2) \}, \, \, \textrm{and}, 
\end{displaymath}
\begin{displaymath}
  (f'_{m-K}(\beta))^{1-(K/m)} \leq f(\beta),
\end{displaymath}
for all $m \geq K$.
The boundedness of the sequence $\{ c_n^{1/n} \}$, where $c_n$ is the total number of animals of size $n$ up to translation, may now be used to show that there exists a constant $L > 1$, which may be chosen uniformly in $\beta \in (0,2(d-1)) - \{ 2(d-2) \}$, such that
$ f'_n(\beta) \leq L^{1/n}f(\beta)$ for all $\beta \in (0,2(d-1)) - \{ 2(d-2) \}$.
The case where $\beta$ assumes the value $2(d-2)$ may be handled by,
for example, altering the way in which the $i$ extra edges are added
to the construction $S$. This case is no more involved than the
general one, and the details are omitted. 
Note now that, for $n,m \in \mathbb{N}$,
\begin{displaymath}
\sigma_{n,m} \leq \sigma'_{n,m} \leq (n + 1)\sigma_{n,m}.
\end{displaymath}
To see the first inequality, note that $\Gamma_{n,m} \subseteq
\Gamma'_{n,m}$; the second, at most $n+1$ translates of a member
$\gamma$ of $\Gamma_n^m$ lie in $\Gamma'_{n,m}$, since some vertex of
$\gamma$ must be mapped to the origin by the translation. We deduce
that $f$ as defined in Theorem \ref{thmone} exists, and satisfies the
bound given in the third part of the Theorem. 

It is easy to show by an induction on $n$, that there is no lattice
animal in $\mathbb{Z}^d$ that has $n$ edges and greater that $2(d-1)n
+ 2d$ outlying edges. Hence $f_n(\beta)$ is equal to zero for $\beta >
2(d-1)$ and $n$ sufficiently large. This implies that $f(\beta)=0$ for
such values of $\beta$. It remains only to prove the fourth part of
Theorem \ref{thmone}.

\subsection{$\bf{f}$ is log-concave}
We complete the proof of Theorem \ref{thmone}, by showing that $f$ is log-concave on $(0,2(d-1))$:
for all $\beta_1,\beta_2 \in (0,2(d-1)), \lambda \in [0.1]$, we have that
\begin{displaymath}
\log f (\lambda \beta_1 + (1-\lambda) \beta_2 ) \geq \lambda \log f (\beta_1)+ (1 - \lambda)\log f (\beta_2).
\end{displaymath}
The technique for proving this result is very much the same as for proving the existence of $f$. We concatenate two large lattice animals of surface-area-to-volume ratio close to $\beta_1$ and $\beta_2$, choosing the relative size of the animals by a weighting determined by $\lambda$. The resulting animal has surface-area-to-volume ratio close to $\lambda \beta_1 + (1- \lambda) \beta_2 $, and the small adjustment required is performed by adding a correction construction. The details of this construction are exactly as in the proof of the existence of $f$. 

So, let $\beta_1, \beta_2, \lambda$ (satisfying the stated conditions) be given, and let $n \in \mathbb{N}$.
Set $n_1 = \lfloor \lambda n \rfloor, n_2 = \lfloor (1 - \lambda) n \rfloor$ and $m_1 = \lfloor \beta_1 \lfloor \lambda n \rfloor \rfloor, m_2 = \lfloor \beta_2 \lfloor(1- \lambda)n \rfloor \rfloor$.
Let $\gamma_1 \in \Gamma'_{n_1,m_1}$ and $\gamma_2 \in \Gamma'_{n_2,m_2}$.
Let $n_3$ and $m_3$ be such that
\begin{displaymath}
\gamma_1 * R_2 * \gamma_2 \in \Gamma'_{n_3,m_3}
\end{displaymath}
Then $n_3 = n_1 + n_2 + 2$ and $m_3 = m_1 + m_2 + 2(d-2)$.
We seek $K$ such that there exists an animal $C$  such that
\begin{displaymath}
\gamma_1 * R_2 * \gamma_2 * R_2 * C \in \Gamma'_{n + K,\lfloor ( \lambda \beta_1 + (1 - \lambda)\beta_2 )(n+K) \rfloor }
\end{displaymath}
The existence of $K$ is established by an almost verbatim rehearsal of the argument proving the existence of $f$.
Let us see how it implies the log-concavity of $f$.
The existence of $K$ and $C$ implies that
\begin{displaymath}
\sigma'_{\lfloor \lambda n \rfloor,\lfloor \beta_1 \lfloor \lambda
n\rfloor \rfloor} \sigma'_{\lfloor (1 - \lambda) n \rfloor,\lfloor
\beta_1 \lfloor (1 - \lambda) n \rfloor \rfloor} \leq
\sigma'_{n+K,\lfloor(\lambda \beta_1 + ( 1 - \lambda) \beta_2)(n+K)\rfloor}  . 
\end{displaymath}
That is,
\begin{displaymath}
 {\Big( f'_{\lfloor \lambda n \rfloor}(\beta_1) \Big)}^{\lfloor \lambda n \rfloor}  {\Big( f'_{\lfloor (1 - \lambda)n \rfloor }(\beta_2) \Big)}^{\lfloor (1 - \lambda)n \rfloor } \leq {\Big( f'_{n+K}(\lambda \beta_1 + (1 - \lambda)\beta_2) \Big) }^{n+K}.
\end{displaymath}
Taking logarithms gives that
\begin{eqnarray}
 & & \frac{ \lfloor \lambda n \rfloor }{n} \log f'_{\lfloor \lambda n \rfloor}(\beta_1) + \frac{\lfloor (1 - \lambda)n \rfloor}{n} \log f'_{\lfloor (1 - \lambda) n \rfloor}(\beta_2) \nonumber \\ 
 & \leq & (1+ (K/n))\log f'_{n+K}(\lambda \beta_1 + (1 - \lambda)\beta_2) . \nonumber
\end{eqnarray}
Taking the limit as $n \to \infty$ gives that
\begin{displaymath}
\lambda \log f(\beta_1) + (1 - \lambda) \log f (\beta_2) \leq \log f( \lambda \beta_1 + (1 - \lambda) \beta_2),
\end{displaymath}
as required.

\section{Acknowledgements} The research in this report was conducted during the period September 2000 to July 2001, during which time, I was supported by a Domus Graduate Scholarship (Competition B) provided by Merton College, Oxford.
I would like to thank the supervisor of my Masters' thesis, Terry Lyons, for hisencouragement and for many stimulating and helpful discussions, and
Colin MacDiarmid and Mathew Penrose for proof-reading an earlier
version of this document. I thank Amir Dembo and John Cardy for their helpful comments, and Neal Madras for providing some valuable suggestions.

\nocite{MR98d:60139}
\nocite{MR88k:60176}
\nocite{MR95m:82076}
\bibliography{pbiblio}

\begin{thebibliography}{10}

\bibitem{MR82b:82048}
Michael Aizenman, Fran{\c{c}}ois Delyon, and Bernard Souillard.
\newblock Lower bounds on the cluster size distribution.
\newblock {\em J. Statist. Phys.}, 23(3):267--280, 1980.

\bibitem{MR86h:82045}
Michael Aizenman and Charles~M. Newman.
\newblock Tree graph inequalities and critical behavior in percolation models.
\newblock {\em J. Statist. Phys.}, 36(1-2):107--143, 1984.

\bibitem{MR1774341}
Rapha{\"e}l Cerf.
\newblock Large deviations for three dimensional supercritical percolation.
\newblock {\em Ast\'erisque}, (267):vi+177, 2000.

\bibitem{DelyonT}
F.~Delyon.
\newblock {\em Taille, forme et nombre des amas dans les problemes de
  percolation}.
\newblock These de 3eme cycle, Universite Pierre et Marie Curie, Paris, 1980.

\bibitem{flesiaetal}
S.~Flesia, D.~S. Gaunt, C.~E. Soteros, and S.~G. Whittington.
\newblock Letter to the editor.
\newblock {\em J. Phys. A}, 25(25):L1169--L1172, 1992.

\bibitem{MR1304212}
S.~Flesia, D.~S. Gaunt, C.~E. Soteros, and S.~G. Whittington.
\newblock Statistics of collapsing lattice animals.
\newblock {\em J. Phys. A}, 27(17):5831--5846, 1994.

\bibitem{grim}
Geoffrey Grimmett.
\newblock {\em Percolation}.
\newblock Springer-Verlag, Berlin, second edition, 1999.

\bibitem{pla}
Alan Hammond.
\newblock {\em Percolation and lattice animals: exponent relations, and
  conditions for $\theta(p_c)=0$}.
\newblock Preprint.

\bibitem{MR98d:60139}
Barry~D. Hughes.
\newblock {\em Random walks and random environments. {V}ol. 2}.
\newblock The Clarendon Press Oxford University Press, New York, 1996.
\newblock Random environments.

\bibitem{ROT}
E.~Janse~van Rensburg, E. J.~Orlandini and M.C. Tesi.
\newblock Collapsing animals.
\newblock {\em J. Phys. A}, 32(32):1567--1584, 1999.

\bibitem{MR58:14852}
Herv{\'e} Kunz and Bernard Souillard.
\newblock Essential singularity in percolation problems and asymptotic behavior
  of cluster size distribution.
\newblock {\em J. Statist. Phys.}, 19(1):77--106, 1978.

\bibitem{MR90c:82051}
N.~Madras, C.~E. Soteros, and S.~G. Whittington.
\newblock Statistics of lattice animals.
\newblock {\em J. Phys. A}, 21(24):4617--4635, 1988.

\bibitem{MR95m:82076}
Neal Madras.
\newblock A rigorous bound on the critical exponent for the number of lattice
  trees, animals, and polygons.
\newblock {\em J. Statist. Phys.}, 78(3-4):681--699, 1995.

\bibitem{MR2001f:82034}
Neal Madras.
\newblock A pattern theorem for lattice clusters.
\newblock {\em Ann. Comb.}, 3(2-4):357--384, 1999.
\newblock On combinatorics and statistical mechanics.

\bibitem{MR89k:82058}
Yigal Meir, Amnon Aharony, and A.~Brooks Harris.
\newblock Percolation in negative field and lattice animals.
\newblock {\em Phys. Rev. B (3)}, 39(1):649--656, 1989.

\bibitem{MR88k:60176}
C.~M. Newman.
\newblock Inequalities for $\gamma$ and related critical exponents in short and
  long range percolation.
\newblock In {\em Percolation theory and ergodic theory of infinite particle
  systems (Minneapolis, Minn., 1984--1985)}, pages 229--244. Springer, New
  York, 1987.

\bibitem{MR1421925}
L.~M. Stratychuk and C.~E. Soteros.
\newblock Statistics of collapsed lattice animals: rigorous results and {M}onte
  {C}arlo simulations.
\newblock {\em J. Phys. A}, 29(22):7067--7087, 1996.

\end{thebibliography}
\bibliographystyle{plain}
\end{document}